\documentclass[a4paper, 11pt]{article}

\usepackage{amsmath}   % From the American Mathematical Society
\usepackage{amssymb}
\usepackage{amsfonts}
\usepackage{wasysym}
\usepackage{dsfont}
\usepackage{graphicx}
\usepackage{subfigure}
\usepackage[latin1]{inputenc}
\newtheorem{thm}{Theorem}
\newtheorem{lem}[thm]{Lemma}
\newtheorem{defn}[thm]{Definition}

\newtheorem{rk}[thm]{Remark}
\newtheorem{cor}[thm]{Corollary}

\begin{document}

%\title{An Alternative Expression for the Joint Eigenvalue Distribution of Uncorrelated Non-central Wishart Matrices with Rank 1 Mean}
\title{Three Problems Related to the Eigenvalues of Complex Non-central Wishart Matrices with a Rank-$1$ Mean}

\author{Prathapasinghe Dharmawansa\\{}\\ \normalsize Department of Statistics, Sequoia Hall, 390 Serra Mall, Stanford University\\ \normalsize Stanford, CA 94305, email: prathapa@stanford.edu}
%Department of Statistics, Sequoia Hall, 390 Serra Mall, Stanford University, Stanford, CA 94305, email: prathapa@stanford.edu
\date{}
\maketitle

\begin{abstract}
Recently, D. Wang \cite{Wang} has devised a new contour integral based method to simplify certain matrix integrals. Capitalizing on that approach, we derive a new expression for the probability density function (p.d.f.) of the joint eigenvalues of a complex non-central Wishart matrix with a rank-$1$ mean. The resulting functional form in turn enables us to use powerful classical orthogonal polynomial techniques in solving three problems related to the non-central Wishart matrix. To be specific, for an $n\times n$ complex non-central Wishart matrix $\mathbf{W}$ with $m$ degrees of freedom ($m\geq n$) and a rank-$1$ mean, we derive  a new expression for the cumulative distribution function (c.d.f.) of the minimum eigenvalue ($\lambda_{\min}$). The c.d.f. is expressed as the determinant of a square matrix, the size of which depends only on the difference $m-n$. This further facilitates the analysis of the microscopic limit for the minimum eigenvalue
which takes the form of the determinant of a square matrix of size $m-n$ with the Bessel kernel. We also develop a moment generating function based approach to derive the p.d.f.  of the random variable $\frac{\text{tr}(\mathbf{W})}{\lambda_{\min}}$, where $\text{tr}(\cdot)$ denotes the trace of a square matrix. This random quantity is of great importance in the so-called smoothed analysis of Demmel condition number. Finally, we find the average  of the reciprocal of the characteristic polynomial $\det[z\mathbf{I}_n+\mathbf{W}],\; |\arg z|<\pi$, where $\mathbf{I}_n$ and $\det[\cdot]$  denote the identity matrix of size $n$ and the determinant, respectively.

\end{abstract}

\section{Introduction}
The eigenvalues of random matrices are known to have far-reaching implications in various scientific disciplines. Finite dimensional properties of the eigenvalues of Wishart type random matrices are of paramount importance in classical multivariate analysis \cite{ Anderson, Robb}, whereas recent multivariate statistical investigations have focused on establishing the asymptotic properties of the eigenvalues \cite{Imj, Alex}. Various links between the eigenvalues of random matrices and statistical physics, combinatorics and integrable systems have been established over the last few decades (see, e.g.,\cite{PJF, Mehta} and references therein). Apart from these areas, random matrices, especially matrices with complex Gaussian elements, have also found new applications in signal processing and wireless communications \cite{Telatar, Verdu}.

The majority of those studies focus on random matrix ensembles derived from zero mean Gaussian matrices. However, random matrices derived from non-zero mean Gaussian matrices have been traditionally an area of interest in multivariate analysis \cite{Anderson, Cons, James, Robb}. Moreover, mathematical objects such as zonal polynomials \cite{Hua, James} and hypergeometric functions of matrix arguments \cite{Herz, Khatri} have been introduced in multivariate analysis literature to facilitate further analysis of such non-central random matrices. Interestingly, these non-central matrices have also been referred to as random matrices with external sources in the literature of physics \cite{Brezin1, Brezin2, Brezin4, Justin}.  In this respect, the classical orthogonal polynomial based characterization of the eigenvalues of random matrices \cite{Mehta} has been further extended to encompass multiple orthogonal polynomials in \cite{Bleher2, Bleher1}. Alternatively, capitalizing on a contour integral approach due to Kazakov \cite{Kazakov}, the authors in \cite{Arous, Brezin1, Brezin2} have introduced a double contour integral representation for the correlation kernel of the eigenvalue point process of non-central random matrices. Some recent contributions on this matter include \cite{ Forrester1, Peter}.

One of the salient features common to those latter studies is that they exclusively focus either on spiked correlation or mean model. It is noteworthy that these two models are mathematically related to each other \cite{Bleher3}. As we are well aware of, the characterization of the joint eigenvalue distribution of non-central random matrices\footnote{Here the term ``non-central random matrices" refers to non-central Gaussian and Wishart matrices.} involves the hypergeometric function of two matrix arguments \cite{James}. It turns out that one of the argument matrices becomes reduced-rank in the presence of a spiked mean/correlation model. Specifically,  when the spike is of rank one, an alternative representation of the hypergeometric function of two matrix arguments has recently been discovered independently by Mo \cite{Mo}, Wang \cite{Wang} and Onatski \cite{Alex}. The key contribution amounts to the representation of the hypergeometric function of two matrix arguments with a rank-$1$ argument matrix in terms of an infinite series involving a single contour integral. This representation has been subsequently used to further characterize the  asymptotic behaviors of the eigenvalues of non-central random matrices \cite{Mo, Wang}.

In this paper, by employing this alternative contour integral representation, we analyze three problems pertaining to the eigenvalues of a finite dimensional complex non-central Wishart matrix with a rank-$1$ mean matrix\footnote{This is also known as the shifted mean chiral Gaussian ensemble with $\beta=2$ (i.e., the complex case) \cite{Peter}.}. Let $0<\lambda_1\leq \lambda_2\leq \cdots\leq \lambda_n$ be the ordered eigenvalues of an $n\times n$ complex non-central Wishart matrix $\mathbf{W}$ with $m$ degrees of freedom and a rank-$1$ mean. We are interested in the following three problems.
\begin{enumerate}
\item The characterization of the cumulative distribution function (c.d.f.) of the minimum eigenvalue of $\mathbf{W}$ as the determinant of a square matrix, the size of which depends on the difference of the degrees of freedom and $n$ (i.e., $m-n$).
\item The statistical characterization of the random quantity $\frac{\text{tr}(\mathbf{W})}{\lambda_1}$ with $\text{tr}(\cdot)$ denoting the trace of a square matrix.
\item The statistical average of the reciprocal of the characteristic polynomial $\det[z\mathbf{I}_n+\mathbf{W}], \;|\arg z|<\pi$, with $\det[\cdot]$ and $\mathbf{I}_n$ denoting the determinant of a square matrix and the $n\times n$ identity matrix, respectively.
\end{enumerate} 

The first problem has a straightforward solution in the form of the determinant of a square matrix of size $n\times n$ \cite{Matthesis}. This stems from the determinant representation of the hypergeometric function of two matrix arguments due to Khatri \cite{Khatri}. However, in certain cases, it is convenient to have an expression with the determinant of a square matrix of size $m-n$ (e.g., when $m=n$). Therefore, in this work, by leveraging the knowledge of classical orthogonal polynomials, we derive an alternative  expression for the c.d.f. of the minimum eigenvalue which involves the determinant of a square matrix of size $m-n+1$. This new form is highly desirable when the difference between $m$ and $n$ is small irrespective of their individual magnitudes. In such a situation, this new expression circumvents the analytical complexities associated with the above straightforward solution which requires to evaluate the determinant of an $n\times n$ square matrix. This key representation, in turn, facilitates the further analysis of the so-called microscopic limit of the minimum eigenvalue (i.e., the limit when $m,n\to\infty$ such that $m-n$ is fixed) which is known to have a determinantal form involving the Bessel kernel \cite{Arous}. 

The random quantity of our interest in the second problem is commonly known as the Demmel condition number in the literature of numerical analysis \cite{Demmel}. As opposed to the case corresponding to the central Wishart matrices \cite{Robb}, $\text{tr}(\mathbf{W})$ and $\frac{\lambda_1}{\text{tr}(\mathbf{W})}$ are no longer statistically independent. Furthermore, a direct Laplace transform relationship between $\frac{\lambda_1}{\text{tr}(\mathbf{W})}$ and the probability density of the minimum eigenvalue of $\mathbf{W}$ does not seem to exist, whereas such a relationship exists among these random quantities in the case of central Wishart matrices \cite{PrathaSIAM, Krish}. Therefore, we introduce a moment generating function (m.g.f.) based framework to solve the second problem. In particular, using a classical orthogonal polynomial approach, we derive the m.g.f.  of the random variable of our interest in terms of a single integral involving the determinant of a square matrix of size $m-n+1$. Upon taking the direct Laplace inversion of the m.g.f. we then obtain an exact expression for the probability desnity function (p.d.f.).  The remarkable fact of having the determinant of a square matrix of size $m-n+1$ makes it suitable to be used when the relative difference between $m$ and $n$ is small. For instance, in the special case of $m=n$, the p.d.f. simplifies to an expression involving a single infinite summation. 

A generalized framework based on the duality between  certain matrix ensembles has been proposed in \cite{Patric} to solve certain problems involving the averages of the reciprocals of characteristic polynomials pertaining to non-central Wishart matrices. However, the third problem of our interest does not seem to be consistent with that framework, since the specific parameters associated with our problem do not satisfy the requirements in \cite{Patric}. Also, it is worth mentioning that this particular problem has not been addressed in a more recent work of Forrester \cite{Forrester1} on the averages of characteristic polynomials for shifted mean chiral Gaussian ensembles. Therefore, again following the classical orthogonal polynomial approach, here we derive a new expression for this particular average. The resultant expression turns out to have a single infinite series. This is not surprising, since in the case of a central Wishart matrix the corresponding answer depends only on the number of characteristic polynomials rather than the size of the random matrix \cite{Boro, Patric, Fodorov}.

The rest of this paper is organized as follows. We begin in Section 2 by deriving a new p.d.f. for the eigenvalues of a complex non-central Wishart matrix with a rank-$1$ mean. In Section $3$ we use the new joint eigenvalue p.d.f. to derive the c.d.f. of the minimum eigenvalue in terms of the determinant of  a square matrix of size $m-n+1$.  Section $4$ addresses the problem of statistical characterization of the random quantity $\frac{\text{tr}(\mathbf{W})}{\lambda_1}$ by deriving corresponding m.g.f. and p.d.f. expressions. Section $5$ is dedicated to deriving the average of the reciprocal of the characteristic polynomial $\det[z\mathbf{I}_n+\mathbf{W}], \;|\arg z|<\pi$.

\section{New Joint Density of the Eigenvalues of a Complex Non-central Wishart Matrix with a Rank-$1$ Mean}
Let us first define the p.d.f. of a complex non-central Wishart matrix.
\begin{defn}
Let $\mathbf{X}\in\mathbb{C}^{m\times n}$ be distributed as $\mathcal{CN}_{n,m}\left(\mathbf{M},\mathbf{I}_m\otimes\mathbf{I}_n\right)$ where $\mathbf{M}\in\mathbb{C}^{n\times n}$ with $m\geq n$. Then $\mathbf{W}=\mathbf{X}^\dagger \mathbf{X}$ has a complex non-central Wishart distribution $\mathcal{W}_n\left(m,\mathbf{I}_n,\mathbf{M}^\dagger\mathbf{M}\right)$ with p.d.f.\footnote{Henceforth, we use $(\cdot)^\dagger$ to denote the conjugate transpose of a matrix.}
\begin{align}
\label{den}
f_{\mathbf{W}}(\mathbf{W})=\frac{e^{-\rm{tr}(\mathbf{M}^\dagger \mathbf{M})}}{\widetilde \Gamma_n(m)}|\mathbf{W}|^{m-n}e^{-\rm{tr}\left(\mathbf{W}\right)}
{}_0\widetilde F_1\left(m;\mathbf{M}^\dagger\mathbf{M}\mathbf{W}\right)
\end{align}
where $\widetilde \Gamma_n(m)=\displaystyle \pi^{\frac{m(m-1)}{2}}\prod_{i=1}^n\Gamma(m-i+1)$ and ${}_0\widetilde F_1\left(\cdot;\cdot\right)$ denotes the complex hypergeometric function of one matrix argument. In particular, for a Hermitian positive definite $n\times n$ matrix $\mathbf{A}$, we have \cite{James}
\begin{equation*}
{}_0\widetilde F_1\left(p;\mathbf{A}\right)=\sum_{k=0}^\infty
\frac{1}{k!}\sum_{\kappa}\frac{C_\kappa(\mathbf{A})}{[p]_\kappa}
\end{equation*}
where $C_\kappa(\cdot)$ is the complex zonal polynomial\footnote{The specific definition of the zonal polynomial is not given here as it is not required in the subsequent analysis. More details of the zonal polynomials can be found in \cite{James, Takemura}.}  which depends through the eigenvalues of the argument matrix $\mathbf{A}$, $\kappa=(k_1,k_2,\ldots,k_n)$, with $k_i$'s being non-negative integers, is a partition of $k$ such that $k_1\geq k_2\geq\ldots\geq k_n\geq 0$ and $\sum_{i=1}^nk_i=k$. Also $[n]_\kappa=\prod_{i=1}^n(n-i+1)_{k_i}$ with $(a)_n=a(a+1)\ldots(a+n-1)$ denoting the Pochhammer symbol.
\end{defn}
The following theorem is due to James \cite{James}.

\begin{thm}
The joint density of the ordered eigenvalues $0<\lambda_1\leq \lambda_2\leq\ldots\leq\lambda_n$ of $\mathbf{W}$ is given by \cite{James}
\begin{align}
\label{eig}
f_{\boldsymbol{\Lambda}}\left(\lambda_1,\lambda_2,\ldots,\lambda_n\right)=K_{m,n}e^{-\rm{tr}\left(\mathbf{M}^\dagger\mathbf{M}\right)}&\Delta_n^2(\boldsymbol{\lambda})\prod_{i=1}^n\lambda_i^{m-n}e^{-\lambda_i} {}_0\widetilde F_1\left(m;\boldsymbol{\Lambda},\mathbf{M}^\dagger\mathbf{M}\right)
\end{align}
where
\begin{equation*}
K_{m,n}=\frac{1}{\prod_{i=1}^n\Gamma(m-i+1)\Gamma(n-i+1)},
\end{equation*}
$\boldsymbol{\Lambda}=\rm{diag}(\boldsymbol{\lambda})$ with $\boldsymbol{\lambda}=\left(\lambda_1,\lambda_2,\ldots,\lambda_n\right)$ and $\Delta_n(\boldsymbol{\lambda})=\prod_{1\leq i<k\leq n}\left(\lambda_k-\lambda_i\right)$. Moreover, ${}_0\widetilde F_1(\cdot;\cdot,\cdot)$ denotes the complex hypergeometric function of two matrix arguments. For Hermitian positive definite $n\times n$ matrices $\mathbf{A}$ and $\mathbf{B}$, we have
\begin{align*}
{}_0\widetilde F_1\left(m;\mathbf{A},\mathbf{B}\right)=\sum_{k=0}^\infty \frac{1}{k!} \sum_{\kappa} \frac{C_\kappa(\mathbf{A})C_\kappa(\mathbf{B})}{[m]_\kappa C_\kappa(\mathbf{I}_n)}.
\end{align*} 
\end{thm}
Now let us focus on simplifying the hypergeometric function in the case of a rank-$1$ mean matrix (i.e., the matrix $\mathbf{M}$ is rank-$1$). To this end, 
 we expand the hypergeometric function term in (\ref{eig}) to yield
\begin{equation}
\label{zonal}
{}_0\widetilde F_1\left(m;\boldsymbol{\Lambda},\mathbf{M}^\dagger\mathbf{M}\right)=
\sum_{k=0}^\infty\frac{1}{k!}\sum_{\kappa}\frac{C_\kappa(\boldsymbol{\Lambda})C_\kappa(\mathbf{M}^\dagger\mathbf{M})}{[m]_\kappa C_\kappa(\mathbf{I}_n)}.
\end{equation}
Since $\mathbf{M}$ is of rank one, clearly the product $\mathbf{M}^\dagger\mathbf{M}$ contains only one non-zero eigenvalue (say $\mu$). This, along with \cite [Corollary 7.2.4]{Robb}, in turn gives that $C_\kappa(\mathbf{M}^\dagger\mathbf{M})=0$, for all partitions of $k$ having more than one non-zero parts. Therefore, only partitions of the form $(k,0,0\cdots,0)$, which we simply denote by $k$, contribute to the summation. In light of this observation, we can simplify (\ref{zonal}) to obtain
\begin{equation}
\label{zonal_simple}
{}_0\widetilde F_1\left(m;\boldsymbol{\Lambda},\mathbf{M}^\dagger\mathbf{M}\right)=
\sum_{k=0}^\infty\frac{1}{(m)_k k!}\left(\prod_{i=0}^{k-1}\frac{1+i}{n+i}\right)C_k(\boldsymbol{\Lambda})\mu^k
\end{equation}
where we have used the fact that $C_k(\mathbf{I}_n)=\prod_{i=0}^{k-1}\frac{n+i}{1+i}$.
Now following Wang \cite{Wang}, we have
\begin{align}
\label{taylor}
\frac{1}{k!}\left(\prod_{i=0}^{k-1}1+i\right)C_k(\boldsymbol{\Lambda})=\frac{1}{2\pi \mathrm{i}}\oint_0\prod_{j=1}^n
\frac{1}{\left(1-z\lambda_j\right)}\frac{{\rm d}z}{z^{k+1}}
\end{align}
where the contour is taken to be a small circle around $0$ with $\frac{1}{\lambda_i} (i=1,2,\ldots,n)$ being exterior of the contour and $\mathrm{i}=\sqrt{-1}$. Substituting (\ref{taylor}) back into (\ref{zonal_simple}) followed by exchanging the summation and integral then gives
%\begin{align*}
%{}_0\widetilde F_1\left(m;\boldsymbol{\Lambda},\mathbf{M}^\dagger\mathbf{M}\right)=\frac{1}{2\pi \mathrm{i}}\oint_0\prod_{j=1}^n
%\frac{1}{\left(1-z\lambda_j\right)}
%\sum_{k=0}^\infty\frac{1}{(m)_k}\left(\prod_{i=0}^{k-1}\frac{1}{n+i}\right)\frac{\mu^k}{z^{k+1}} {\rm d}z,
%\end{align*} ,
\begin{align*}
{}_0\widetilde F_1\left(m;\boldsymbol{\Lambda},\mathbf{M}^\dagger\mathbf{M}\right)=\frac{1}{2\pi \mathrm{i}}\oint_0\prod_{j=1}^n
\frac{1}{\left(1-z\lambda_j\right)}
\sum_{k=0}^\infty\frac{1}{(m)_k(n)_k}\frac{\mu^k}{z^{k+1}} {\rm d}z
\end{align*} 
where we have used the relation $\prod_{i=0}^{k-1}(n+i)=(n)_k$.
Since there exists an integer $N$ such that $n=N+1$, we can rewrite  the above equation as
 \begin{align*}
{}_0\widetilde F_1\left(m;\boldsymbol{\Lambda},\mathbf{M}^\dagger\mathbf{M}\right)
%\frac{1}{2\pi \mathrm{i}}\oint_0\prod_{j=1}^n
%\frac{1}{\left(1-z\lambda_j\right)}
%\sum_{k=0}^\infty\frac{1}{(m)_k(N+1)_k}\frac{\mu^k}{z^{k+1}} {\rm d}z\nonumber\\
%&=\frac{N!(m-1)!}{2\pi \mathrm{i}}\oint_0\prod_{j=1}^n
%\frac{1}{\left(1-z\lambda_j\right)}
%\sum_{k=0}^\infty\frac{1}{\Gamma(m+k)\Gamma(N+k+1)}\frac{\mu^k}{z^{k+1}} {\rm d}z\nonumber\\
& =\frac{N!(m-1)!}{2\pi \mathrm{i}}\oint_0\prod_{j=1}^n
\frac{1}{\left(1-z\lambda_j\right)}
\sum_{k=N}^\infty\frac{1}{\Gamma(m+k-N)k!}\frac{\mu^{k-N}}{z^{k-N+1}} {\rm d}z\nonumber\\
&=\frac{N!(m-1)!}{(m-n)!}\frac{1}{2\pi \mathrm{i}}\oint_0\prod_{j=1}^n
\frac{1}{\left(1-z\lambda_j\right)}\left\{
\sum_{k=0}^\infty\frac{1}{k!(m-N)_k}\frac{\mu^{k-N}}{z^{k-N+1}} \right.\nonumber\\
& \hspace{5.5cm} -\left.\sum_{k=0}^{N-1}\frac{1}{k!(m-N)_k}\frac{\mu^{k-N}}{z^{k-N+1}} \right\}{\rm d}z.
\end{align*}
Clearly, the second summation evaluates to zero, since the integrand is an analytic function. Therefore, we can further simplify the above equation to yield
\begin{align*}
{}_0\widetilde F_1\left(m;\boldsymbol{\Lambda},\mathbf{M}^\dagger\mathbf{M}\right)
%\frac{N!(m-1)!}{(m-n)!}\frac{1}{2\pi \mathrm{i}}\oint_0\prod_{j=1}^n
%\frac{1}{\left(1-z\lambda_j\right)}
%\sum_{k=0}^\infty\frac{1}{k!(m-N)_k}\frac{\mu^{k-N}}{z^{k-N+1}} {\rm d}z\nonumber\\
%&=\frac{N!(m-1)!}{(m-n)!\;\mu^N}\frac{1}{2\pi \mathrm{i}}\oint_0\prod_{j=1}^n
%\frac{1}{\left(1-z\lambda_j\right)}z^{N-1}
%\sum_{k=0}^\infty\frac{1}{k!(m-N)_k}\frac{\mu^{k}}{z^{k}} {\rm d}z\nonumber\\
&=\frac{N!(m-1)!}{(m-n)!\;\mu^N}\frac{1}{2\pi \mathrm{i}}\oint_0\prod_{j=1}^n
\frac{1}{\left(1-z\lambda_j\right)}z^{N-1}
{}_0F_1\left(m-N;\frac{\mu}{z}\right){\rm d}z,
\end{align*}
from which we obtain after the change of variable 
\begin{align}
\label{hyp_simp}
{}_0\widetilde F_1\left(m;\boldsymbol{\Lambda},\mathbf{M}^\dagger\mathbf{M}\right)=\frac{N!(m-1)!}{(m-n)!\;\mu^N}\frac{1}{2\pi \mathrm{i}}\oint_\infty\prod_{j=1}^n
\frac{1}{\left(s-\lambda_j\right)}
{}_0F_1\left(m-N;\mu s\right){\rm d}s
\end{align}
where all the $\lambda_i$'s except $0$ lie inside the contour. Finally, using (\ref{hyp_simp}) in (\ref{eig}) gives the new contour integral representation of the joint p.d.f. of $\lambda_1,\lambda_2,\ldots,\lambda_n$ as
\begin{align*}
f_{\boldsymbol{\Lambda}}\left(\lambda_1,\lambda_2,\ldots,\lambda_n\right)
=K_{m,n}\frac{N!(m-1)!}{(m-n)!}
\frac{e^{-\mu}}{\mu^N}
\frac{1}{2\pi \mathrm{i}}\oint_\infty
{}_0F_1\left(m-N;\mu s\right)
\prod_{j=1}^n
\frac{\lambda_j^{m-n}}{\left(s-\lambda_j\right)}e^{-\lambda_j}
\Delta_n^2(\boldsymbol{\lambda}){\rm d}s.
\end{align*}
One can evaluate the above contour integral to obtain the following new representation for the distribution of the eigenvalues of $\mathbf{W}$. 
\begin{cor}
Let $\mathbf{W}\sim \mathcal{W}_n\left(m,\mathbf{I}_n,\mathbf{M}^\dagger\mathbf{M}\right)$, where $\mathbf{M}$ is rank-$1$ and $\mathrm{tr}(\mathbf{M}^\dagger\mathbf{M})=\mu$. Then the joint density of the eigenvalues of $\mathbf{W}$ is given by
\begin{align}
\label{newden}
&f_{\boldsymbol{\Lambda}}\left(\lambda_1,\lambda_2,\ldots,\lambda_n\right)=\mathcal{K}_{n,\alpha}
\frac{e^{-\mu}}{\mu^{n-1}}
\prod_{i=1}^n\lambda_i^{\alpha}e^{-\lambda_i}
\Delta_n^2(\boldsymbol{\lambda})
\sum_{k=1}^n\frac{{}_0F_1\left(\alpha+1;\mu \lambda_k\right)}{\displaystyle\prod_{\substack{i=1\\i\neq k}}^n\left(\lambda_k-\lambda_i\right)}
\end{align}
where
\begin{equation*}
\mathcal{K}_{n,\alpha}=K_{n+\alpha,n}\frac{(n-1)!(n+\alpha-1)!}{\alpha!}
\end{equation*}
with $\alpha=m-n$.
\end{cor}
Let us now see how to derive a new expression for the c.d.f of the minimum eigenvalue of a complex non-central Wishart matrix with a rank-$1$ mean starting from the joint p.d.f. given above.

Before proceeding, it is worth mentioning the following preliminary results and definitions.
\begin{defn}
For $\rho>-1$, the generalized Laguerre polynomial of degree $M$, $L^{(\rho)}_M(z)$, is given by \cite{Szego}
\begin{equation}
\label{lagdef}
L^{(\rho)}_M(z)=\frac{(\rho+1)_M}{M!}\sum_{j=0}^{M}\frac{(-M)_j}{(\rho+1)_j}\frac{z^j}{j!},
\end{equation}
with the kth derivative satisfying
\begin{equation}
\label{lagderi}
\frac{d^k}{dz^k}L^{(\rho)}_M(z)=(-1)^kL^{(\rho+k)}_{M-k}(z).
\end{equation}
Also $L^{(\rho)}_M(z)$ satisfies the following contiguous relationship
\begin{equation}
\label{contg}
L^{(\rho-1)}_M(z)=L^{(\rho)}_M(z)-L^{(\rho)}_{M-1}(z).
\end{equation}
\end{defn}

\begin{defn}
For a negative integer $-M$, we have the following relation
\begin{equation}
\label{poch}
(-M)_j=\left\{\begin{array}{cc}
(-1)^j\frac{M!}{(M-j)!}& \text{for $j\leq M$}\\
0 & \text{for $j>M$}.
\end{array}\right.
\end{equation}
\end{defn}

\begin{lem}\label{lag}
Following \cite[Eq. 7.414.7]{Grad} and \cite[Corollary 2.2.3]{Askey}, for $j,k\in\{0,1,2,\cdots\}$, we can establish
\begin{equation*}
\int_0^\infty x^j e^{-x} L^{(k)}_M(x){\rm d} x=\frac{j!}{M!}(k-j)_M.
\end{equation*}
\end{lem}

The following compact notation has been used to represent the determinant of an $M\times M$ block matrix:
\begin{equation*}
\det\left[a_{i,1}\;\;\; b_{i,j-1}\right]_{\substack{i=1,2,\cdots,M\\
j=2,3,\cdots,M}}=\left|\begin{array}{ccccc}
a_{1,1} & b_{1,1} & b_{1,2}& \cdots & b_{1,M-1}\\
a_{2,1} & b_{2,1} & b_{2,2}& \cdots & b_{2,M-1}\\
\vdots & \vdots & \vdots & \ddots & \vdots\\
a_{M,1} & b_{M,1} & b_{M,2}& \cdots & b_{M,M-1}
\end{array}
\right|.
\end{equation*}

\section{Cumulative Distribution of the Minimum Eigenvalue}
Here we derive a new expression for the c.d.f. of the minimum eigenvalue $\lambda_{\min}$ of $\mathbf{W}$ with a rank-$1$ mean.

By definition, the c.d.f. of $\lambda_{\min}$ is given by
\begin{equation}
\label{cdf}
F_{\lambda_{\min}}(x)=\Pr\left(\lambda_1<x\right)=1-\Pr\left(\lambda_1\geq x\right)
\end{equation}
where
\begin{equation}
\label{multi_integral_one}
\Pr\left(\lambda_1\geq x\right)=\int_{x\leq \lambda_1\leq \lambda_2\leq \cdots\leq \lambda_n<\infty}
f_{\boldsymbol{\Lambda}}\left(\lambda_1,\lambda_2,\ldots,\lambda_n\right)  {\rm d}\lambda_1{\rm d}\lambda_2\cdots{\rm d}\lambda_n.
\end{equation}
The following theorem gives the c.d.f. of $\lambda_{\min}$.
\begin{thm}
Let $\mathbf{W}\sim \mathcal{W}_n\left(m,\mathbf{I}_n,\mathbf{M}^\dagger\mathbf{M}\right)$, where $\mathbf{M}$ is rank-$1$ and $\mathrm{tr}(\mathbf{M}^\dagger\mathbf{M})=\mu$. Then the c.d.f. of the minimum eigenvalue of $\mathbf{W}$ is given by
\begin{align}
\label{cdffinal}
F_{\lambda_{\min}}(x)& =1-(n+\alpha-1)!\;e^{-nx}\det\left[(-\mu)^{i-1} \psi_i(\mu, x)
\;\;\; L_{n+i-j}^{(j-2)}(-x)\right]_{\substack{i=1,2,\cdots,\alpha+1\\j=2,3,\cdots,\alpha+1}}
\end{align}
where $\alpha=m-n$,
\begin{align*}
\psi_i(\mu,x)=\frac{1}{(\alpha+i+n-2)!}\sum_{k=0}^\infty \frac{(x\mu)^k{}_1F_1\left(\alpha+k;\alpha+n+i+k-1;-\mu\right)}{k!(\alpha+i+n-1)_k},
\end{align*}
and ${}_1F_1(a;c;z)$ is the confluent hypergeometric function of the first kind.
\end{thm}

%\begin{rk}
%Alternatively, we can express $\psi_i(\mu,x)$ as
%\begin{align}
%\psi_i(\mu,x)=\frac{e^{-\mu}}{(\alpha+i+n-2)!}\Phi_3\left(n+i-1,n+\alpha+i-1;\mu,x\mu\right)
%\end{align}
%\end{rk}

{\bf{Proof:}}
Since the joint p.d.f. is symmetric in $\lambda_1,\lambda_2,\cdots,\lambda_n$, we can write (\ref{multi_integral_one}) as
\begin{equation*}
\Pr\left(\lambda_1\geq x\right)=\frac{1}{n!}\int_{[x,\infty)^n}f_{\boldsymbol{\Lambda}}\left(\lambda_1,\lambda_2,\ldots,\lambda_n\right)  {\rm d}\lambda_1{\rm d}\lambda_2\cdots{\rm d}\lambda_n
\end{equation*}
from which we obtain upon using the variable transformations
\begin{equation*}
\Pr\left(\lambda_1\geq x\right)=\frac{1}{n!}\int_{[0,\infty)^n}f_{\boldsymbol{\Lambda}}\left(\lambda_1+x,\lambda_2+x,\ldots,\lambda_n+x\right)  {\rm d}\lambda_1{\rm d}\lambda_2\cdots{\rm d}\lambda_n.
\end{equation*}
Now it is convenient to use (\ref{newden}) to arrive at
\begin{align}
\label{eqsum}
\Pr\left(\lambda_1\geq x\right)= \frac{\mathcal{K}_{n,\alpha}}{n!}\frac{e^{-\mu-nx}}{\mu^{n-1}}
 \sum_{k=1}^n \int_{[0,\infty)^n}& 
\frac{{}_0F_1\left(\alpha+1;\mu (\lambda_k+x)\right)}{\prod_{\substack{i=1\\i\neq k}}^n\left(\lambda_k-\lambda_i\right)} \nonumber\\
&  \qquad \times
\prod_{i=1}^n(\lambda_i+x)^{\alpha}e^{-\lambda_i}
\Delta_n^2(\boldsymbol{\lambda}) {\rm d}\lambda_1{\rm d}\lambda_2\cdots{\rm d}\lambda_n.
\end{align}
Since each term in the above summation contributes the same amount, we may write (\ref{eqsum}) as
\begin{align*}
\Pr\left(\lambda_1\geq x\right)= \frac{\mathcal{K}_{n,\alpha}}{(n-1)!}\frac{e^{-\mu-nx}}{\mu^{n-1}}
 \int_{[0,\infty)^n} &
\frac{{}_0F_1\left(\alpha+1;\mu (\lambda_1+x)\right)}{\prod_{i=2}^n\left(\lambda_1-\lambda_i\right)} \nonumber\\
& \qquad \qquad\times
\prod_{i=1}^n(\lambda_i+x)^{\alpha}e^{-\lambda_i}
\Delta_n^2(\boldsymbol{\lambda}) {\rm d}\lambda_1{\rm d}\lambda_2\cdots{\rm d}\lambda_n.
\end{align*}
Noting the fact that
\begin{equation*}
\Delta_n^2(\boldsymbol{\lambda})=\prod_{i=2}^n(\lambda_1-\lambda_i)^2 \Delta_{n-1}^2(\boldsymbol{\lambda})
\end{equation*}
with $\Delta_{n-1}^2(\boldsymbol{\lambda})=\prod_{2\leq k<j\leq n}(\lambda_j-\lambda_k)^2$, we can rewrite the above multiple integral after some algebraic manipulation as
\begin{align*}
\Pr\left(\lambda_1\geq x\right)=& \frac{\mathcal{K}_{n,\alpha}}{(n-1)!}\frac{e^{-\mu-nx}}{\mu^{n-1}}
 \int_{[0,\infty)}
{}_0F_1\left(\alpha+1;\mu (\lambda_1+x)\right) (\lambda_1+x)^\alpha e^{-\lambda_1}\nonumber\\
& \qquad \times
\left(\int_{[0,\infty)^{n-1}}\prod_{i=2}^n(\lambda_1-\lambda_i)(\lambda_i+x)^{\alpha}e^{-\lambda_i}
\Delta_{n-1}^2(\boldsymbol{\lambda}) {\rm d}\lambda_2\cdots{\rm d}\lambda_n\right) \; {\rm d}\lambda_1.
\end{align*}
Now it is convenient to relabel the variables as $\lambda_1=\lambda$ and $\lambda_i=y_{i-1}, i=2,3,\cdots,n$, to obtain
%\begin{align}
%\Pr\left(\lambda_1\geq x\right)=& \frac{\mathcal{K}_{n,\alpha}}{(n-1)!}\frac{e^{-\mu-nx}}{\mu^{n-1}}
% \int_{[0,\infty)}
%{}_0F_1\left(\alpha+1;\mu (\lambda+x)\right) (\lambda+x)^\alpha e^{-\lambda}\nonumber\\
%& \qquad \times
%\left(\int_{[0,\infty)^{n-1}}\prod_{i=1}^{n-1}(\lambda-y_i)(y_i+x)^{\alpha}e^{-y_i}
%\Delta_{n-1}^2(\boldsymbol{y}) {\rm d}y_2\cdots{\rm d}y_n\right) \; {\rm d}\lambda
%\end{align}
%which we can rewrite as
\begin{align}
\label{cdfin}
\Pr\left(\lambda_1\geq x\right)= \frac{\mathcal{K}_{n,\alpha}}{(n-1)!}\frac{e^{-\mu-nx}}{\mu^{n-1}}
 \int_{[0,\infty)}&
{}_0F_1\left(\alpha+1;\mu (\lambda+x)\right) (\lambda+x)^\alpha e^{-\lambda}
\nonumber\\
& \qquad \qquad \quad \times
 (-1)^{(n-1)\alpha}Q_{n-1}\left(\lambda,-x,\alpha\right) {\rm d}\lambda
\end{align}
where
\begin{equation}
\label{Qintdef}
Q_{n}\left(a,b,\alpha\right):=\int_{[0,\infty)^n}\prod_{i=1}^{n}(a-y_i)(b-y_i)^{\alpha}e^{-y_i}
\Delta_{n}^2(\boldsymbol{y}) {\rm d}y_1{\rm d} y_2\cdots{\rm d}y_n.
\end{equation}
As shown in the Appendix, we can solve the above multiple integral in closed form giving
\begin{align}
\label{q1}
Q_{n}\left(a,b,\alpha\right) =\frac{\overline{\mathcal{K}}_{n,\alpha}}{(b-a)^\alpha} \det\left[L_{n+i-1}^{(0)}(a)\;\;\; L_{n+i+1-j}^{(j-2)}(b)\right]_{\substack{i=1,2,\cdots,\alpha+1\\j=2,3,\cdots,\alpha+1}}
\end{align}
where
\begin{equation*}
\overline{\mathcal{K}}_{n,\alpha}=(-1)^{n+\alpha(n+\alpha)}\frac{\prod_{i=1}^{\alpha+1}(n+i-1)!\prod_{i=0}^{n-1}i!(i+1)!}{\prod_{i=1}^{\alpha-1}i!}.
\end{equation*}
Therefore, using (\ref{q1}) in (\ref{cdfin}) with some algebraic manipulation, we have
\begin{align*}
\Pr\left(\lambda_1\geq x\right)=(-1)^{n+1}\frac{(n+\alpha-1)!}{\alpha!}\frac{e^{-\mu-nx}}{\mu^{n-1}}&\int_0^\infty
{}_0F_1\left(\alpha+1;\mu (\lambda+x)\right)e^{-\lambda}\nonumber\\
& \times \det\left[ L_{n+i-2}^{(0)}(\lambda)\;\;\; L_{n+i-j}^{(j-2)}(-x)\right]_{\substack{i=1,2,\cdots,\alpha+1\\j=2,3,\cdots,\alpha+1}} {\rm d}\lambda.
\end{align*}
Observing that only the first column of the determinant contains the variable $\lambda$, we can further simplify the above integral to yield
\begin{align}
\label{mineigsplit}
\Pr\left(\lambda_1\geq x\right)&=(-1)^{n+1}\frac{(n+\alpha-1)!}{\alpha!}\frac{e^{-\mu-nx}}{\mu^{n-1}} \det\left[\zeta_i(x)\;\;\; L_{n+i-j}^{(j-2)}(-x)\right]_{\substack{i=1,2,\cdots,\alpha+1\\j=2,3,\cdots,\alpha+1}} 
\end{align}
where
\begin{align*}
\zeta_i(x)=\int_0^\infty
{}_0F_1\left(\alpha+1;\mu (\lambda+x)\right)e^{-\lambda} L_{n+i-2}^{(0)}(\lambda){\rm d}\lambda.
\end{align*}
The remaining task is to evaluate the above integral, which does not seem to have a simple closed-form solution. Therefore, we expand the hypergeometric function with its equivalent power series and use some algebraic manipulation to arrive at
\begin{align*}
\zeta_i(x)=\sum_{l=0}^\infty\sum_{k=0}^\infty \frac{\mu^{l+k} x^{k}}{k!l!(\alpha+1)_{l+k}}\int_0^\infty
 \lambda^l e^{-\lambda} L_{n+i-2}^{(0)}(\lambda){\rm d}\lambda.
\end{align*}
This integral can be solved with the help of Corollary \ref{lag} to yield
\begin{align*}
\zeta_i(x)=\sum_{l=0}^\infty\sum_{k=0}^\infty \frac{\mu^{l+k} x^{k}}{k!(\alpha+1)_{l+k}}\frac{(-l)_{n+i-2}}{(1)_{n+i-2}}.
\end{align*}
Following (\ref{poch}), we observe that the quantity $(-l)_{n+i-2}$ is non-zero only when $l\geq n+i-2$. Therefore, we shift the summation index $l$ with some algebraic manipulation to yield
\begin{align}
\label{zetaans}
\zeta_i(x)=\frac{(-1)^{n+i}\mu^{n+i-2}\alpha !}{(n+i+\alpha-2)!}
\sum_{k=0}^\infty
\frac{(x\mu)^k}{k! (n+i+\alpha-1)_k} {}_1F_1(n+i-1;n+\alpha+i+k-1;\mu)
\end{align}
where we have used the relation
\begin{align*}
(\alpha+1)_{k+i+n+l-2}=\frac{(\alpha+i+n-2)!}{\alpha!}(\alpha+i+n+k-1)_l (\alpha+i+n-1)_k.
\end{align*}
Substituting (\ref{zetaans}) back into (\ref{mineigsplit}) with some algebra then gives 
\begin{align}
\label{hypeq}
\Pr\left(\lambda_1\geq x\right)
& = (n+\alpha-1)!e^{-nx}\det\left[(-\mu)^{i-1}\psi_i(\mu,x)
\;\;\; L_{n+i-j}^{(j-2)}(-x)\right]_{\substack{i=1,2,\cdots,\alpha+1\\j=2,3,\cdots,\alpha+1}}
\end{align}
where we have used the Kummer relation \cite{Erdelyi}
\begin{equation*}
{}_1F_1(a;c;z)=e^z {}_1F_1(c-a;c,-z).
\end{equation*}
Finally, using (\ref{hypeq}) in (\ref{cdf}) gives the c.d.f. of the minimum eigenvalue which concludes the proof.
\begin{rk}
Alternatively, we can express $\psi_i(\mu, x)$ as
\begin{align}
\psi_i(\mu,x)=\frac{e^{-\mu}}{(\alpha+i+n-2)!}\Phi_3\left(n+i-1,n+\alpha+i-1;\mu,x\mu\right)
\end{align}
where
\begin{equation*}
\Phi_3(a,c;x,y)=\sum_{i=0}^\infty\sum_{j=0}^\infty \frac{(a)_i}{(c)_{i+j} i!j!} x^iy^j
\end{equation*}
is the confluent hypergeometric function of two variables \cite[Eq. 5.7.1.23]{Erdelyi}.
\end{rk}

In the special case of $\alpha=0$ (i.e., $m=n$), (\ref{cdffinal}) admits the following simple form
\begin{align}
\label{cdfjmva}
F_{\lambda_{\min}}(x)& =1-e^{-nx}\sum_{k=0}^\infty \frac{(x\mu)^k}{k!(n)_k}{}_1F_1(k;n+k;-\mu)\nonumber\\
&=1-e^{-\mu-nx}\Phi_3\left(n,n;\mu,x\mu\right)
\end{align}
which coincides with what we have derived in \cite[Eq. 32/39]{Pratha} purely based on a matrix integral approach.\footnote{Since the results given in \cite{Pratha} are valid for an arbitrary covariance matrix with $\alpha=0$, one has to assume the identity covariance to obtain the above results.}

In addition, it is not difficult to show that, for $\mu=0$, (\ref{cdffinal}) simplifies to
\begin{equation}
F_{\lambda_{\min}}(x) =1-e^{-nx}\det\left[L_{n+i-j}^{(j-1)}(-x)\right]_{i,j=1,2,\cdots,\alpha}.
\end{equation}

\begin{rk}
Although we can show that the microscopic limit of the above expression (\ref{cdffinal}) for the c.d.f of the minimum eigenvalue takes the form of a determinant of size $\alpha$ involving the Bessel kernel, we omit a detailed analysis here as this particular asymptotic limit is well known in the literature \cite{Arous}.
\end{rk}

Having analyzed the behavior of the minimum eigenvalue of $\mathbf{W}$, let us now move on to determine the distribution of the random variable  $\frac{\rm{tr}(\mathbf{W})}{\lambda_1}$.

\section{The Distribution of $\frac{\rm{tr}(\mathbf{W})}{\lambda_1}$}
Here we study the distribution of the quantity
\begin{equation}
\label{vint}
V=\frac{\rm{tr}(\mathbf{W})}{\lambda_1}=\frac{\sum_{j=1}^n\lambda_j}{\lambda_1}.
\end{equation}
It turns out that this quantity is intimately related to the distribution of the minimum eigenvalue of $\mathbf{W}$ given the constraint $\rm{tr}(\mathbf{W})=1$ (i.e., fixed trace) \cite{Profchen}. To be precise, the latter is distributed as $\frac{1}{V}$. Apart from that, the most notable application of the distribution of $V$ is the so-called ``smoothed analysis of condition numbers" \cite{Spielman}. For a given function $g:\mathbb{C}^{m\times n}\to \mathbb{R}_+$ (e.g., the $2-$norm condition number), $\mathbf{A}\sim \mathcal{CN}_{m,n}(\mathbf{M}, \sigma^2\mathbf{I}_m\otimes \mathbf{I}_n)$ with $0<\sigma\leq 1$ and $\mathbf{M}\in\mathbb{C}^{m\times n}$ being arbitrary such that either $\rm{tr}\left(\mathbf{M}^\dagger\mathbf{M}\right)=1$ or $||\mathbf{M}||_2\leq \sqrt{n}$ is satisfied,  under the smoothed analysis framework, a typical problem is to study the behavior of \cite{Mario, Sankar, Cucker, Burg}
\begin{equation}
\label{smooth}
\sup_{\mathbf{M}}\mathrm{E}_{\mathbf{A}}\left(g(\mathbf{A})\right)
\end{equation}
where $\rm{E}_{\mathbf{A}}(\cdot)$ and $||\cdot||_2$ denote the mathematical expectation with respect to $\mathbf{A}$ and the $2$-norm, respectively. For mathematical tractability, sometimes it is assumed that the matrix $\mathbf{M}$ is of rank one \cite{Mario}. Bounds on the quantity (\ref{smooth}) have been derived in the literature when $g(\mathbf{A})$ defines various condition numbers (see, e.g., \cite{Sankar, Cucker, Burg} and references therein). Among those condition numbers, the one introduced by James Demmel \cite{Demmel} plays an important role in understanding the behaviors of other condition numbers arising in different contexts. For a rectangular matrix $\mathbf{X}\in\mathbb{C}^{m\times n}$, the function $g$ defined by \cite{Cucker}
\begin{equation}
\label{dem_mat}
g(\mathbf{X})=||\mathbf{X}||_F||\mathbf{X}^*||_2
\end{equation}
with $||\cdot||_F$ denoting the Frobenius norm and $\mathbf{X}^*$ denoting the Moore-Penrose inverse, gives the Demmel condition number\footnote{This is the extension of the condition number definition given in \cite{Demmel} to  $m\times n$ rectangular matrices.}. 
In particular, for the matrix of our interest $\mathbf{X}\sim \mathcal{CN}_{m,n}(\mathbf{M},\mathbf{I}_m\otimes \mathbf{I}_n)$ with $m\geq n$, (\ref{dem_mat}) specializes to
\begin{equation*}
g(\mathbf{X})=\frac{\sum_{j=1}^n\lambda_j}{\lambda_1}=V
\end{equation*}
where $\lambda_1\leq \lambda_2\leq\cdots\leq \lambda_n$ are the ordered eigenvalues of $\mathbf{W}=\mathbf{X}^\dagger\mathbf{X}$. In light of  these developments we can clearly see that the distribution of $V$ is of great importance in performing the smoothed analysis on the Demmel condition number.

Having understood the importance of the variable $V$ in (\ref{vint}), we now focus on deriving its p.d.f when the matrix $\mathbf{W}$ has a rank-$1$ mean. For this purpose, here we adopt an approach based on the m.g.f. of $V$. We have the following key result.

\begin{thm}
Let $\mathbf{W}\sim\mathcal{W}_n\left(m,\mathbf{I}_n,\mathbf{M}^\dagger\mathbf{M}\right)$, where $\mathbf{M}$ is rank-$1$ and $\mathrm{tr}(\mathbf{M}^\dagger\mathbf{M})=\mu$. Then the p.d.f. of $V$ is given by
\begin{align}
\label{pdfexact}
 f^{(\alpha)}_V(v) =(n-1)!
\frac{e^{-\mu}}{v^{n(n+\alpha)}} &
\mathcal{L}^{-1}\Biggl\{
\frac{e^{-ns}}{s^{(n-1)(n+\alpha+1)}}
\det\left[\left(-\frac{\mu}{ sv}\right)^{i-1} 
\phi_i(\mu,s,v)
\;\;\; L^{(j)}_{n+i-1-j}(-s)
\right]_{\substack{i=1,2,\cdots,\alpha+1\\
j=2,3,\cdots,\alpha+1}}
\Biggr\}
\end{align}
where
\begin{align*}
\phi_i(\mu,s,v)&=
\sum_{k=0}^\infty \frac{a_i(k)}{k!}
\left(\frac{\mu}{s v}\right)^k {}_1F_1\left(n^2+n\alpha+k+i-1;n+i+k+\alpha-1;\frac{\mu}{v}\right)\\
%\xi_i(k,z)&={}_1F_1\left(n^2+n\alpha+k+i-1;n+i+k+\alpha-1;z\right)\\
a_i(k)&=(n+i-1)\frac{(n^2+n\alpha+i-2)!}{(n+i+\alpha-2)!}\frac{(n+i)_k(n+i-2)_k (n^2+n\alpha+i-1)_k}{(n+i-1)_k(n+i+\alpha-1)_k}
\end{align*}
and $\mathcal{L}^{-1}(\cdot)$ denotes the inverse Laplace transform.
\end{thm}
{\bf{Proof:}} By definition, the m.g.f. of $V$ can be written as
\begin{align*}
\mathfrak{M}_V(s)={\rm{E}}_{\boldsymbol{\Lambda}}\left(e^{-s\frac{\sum_{j=1}^n\lambda_j}{\lambda_1}}\right),\;\;\; \Re(s)\geq 0,
\end{align*}
which has the following multiple integral representation
\begin{align*}
\mathfrak{M}_V(s)=e^{-s}\int_{0\leq \lambda_1\leq\lambda_2\leq\cdots\leq \lambda_n<\infty} e^{-s\frac{\sum_{j=2}^n\lambda_j}
{\lambda_1}}
f_{\boldsymbol{\Lambda}}(\lambda_1,\lambda_2,\cdots,\lambda_n)
{\rm d}\lambda_1{\rm d}\lambda_2\cdots {\rm d}\lambda_n.
\end{align*}
Since the argument of the exponential function is symmetric in $\lambda_2,\cdots,\lambda_n$, it is convenient to introduce the substitution $\lambda_1=x$ and rewrite the multiple integral, keeping the integration with respect to $x$ last, as
\begin{align}
\label{int_split}
\mathfrak{M}_V(s)=e^{-s}\int_0^\infty \int_{x\leq\lambda_2\leq\cdots\leq \lambda_n<\infty} e^{-s\frac{\sum_{j=2}^n\lambda_j}
{x}}
f_{\boldsymbol{\Lambda}}(x,\lambda_2,\cdots,\lambda_n)
{\rm d}\lambda_2\cdots {\rm d}\lambda_n{\rm d}x.
\end{align}
To be consistent with the above setting, we may restructure the joint p.d.f. of $\boldsymbol{\Lambda}$ given in (\ref{newden}) as
\begin{align}
\label{res_den}
f_{\boldsymbol{\Lambda}}(x,\lambda_2,\cdots,\lambda_n)=\mathcal{K}_{n,\alpha}
\frac{e^{-\mu}}{\mu^{n-1}} x^\alpha e^{-x}& 
\prod_{i=2}^n\lambda_i^{\alpha}e^{-\lambda_i}(x-\lambda_i)^2
\Delta_{n-1}^2(\boldsymbol{\lambda})\nonumber\\
&\hspace{-4mm} \times 
\left(\frac{{}_0F_1\left(\alpha+1;\mu x\right)}{\displaystyle\prod_{i=2}^n\left(x-\lambda_i\right)}+\sum_{k=2}^n\frac{{}_0F_1\left(\alpha+1;\mu \lambda_k\right)}{(\lambda_k-x)\displaystyle\prod_{\substack{i=2\\i\neq k}}^n\left(\lambda_k-\lambda_i\right)}\right)
\end{align}
where we have used the decomposition $\Delta_n^2(\boldsymbol{\lambda})=(x-\lambda_i)^2
\Delta_{n-1}^2(\boldsymbol{\lambda})$. Now we use (\ref{res_den}) in (\ref{int_split}) with some algebraic manipulation to obtain
\begin{equation}
\label{mgffact}
\mathfrak{M}_V(s)=\mathfrak{P}(s)+\mathfrak{S}(s)
\end{equation}
where
\begin{align}
\label{Pdef}
\mathfrak{P}(s)&=\mathcal{K}_{n,\alpha}
\frac{e^{-\mu-s}}{\mu^{n-1}}\int_0^\infty
e^{-x}x^\alpha{}_0F_1\left(\alpha+1;\mu x\right)\nonumber\\
& \qquad \quad\times \left( \int_{x\leq \lambda_2\leq\cdots\leq \lambda_n<\infty}
\prod_{i=2}^n e^{-\left(1+\frac{s}{x}\right)\lambda_i}\lambda_i^\alpha (x-\lambda_i)\Delta_{n-1}^2(\boldsymbol{\lambda})
{\rm d}\lambda_2\cdots{\rm d}\lambda_n\right)\; {\rm d}x
\end{align}
and
\begin{align}
\label{Sdef}
\mathfrak{S}(s)=\mathcal{K}_{n,\alpha}
\frac{e^{-\mu-s}}{\mu^{n-1}}
\int_0^\infty e^{-x}x^{\alpha}& \Biggl(\int_{x\leq \lambda_2\leq\cdots\leq \lambda_n<\infty}
\sum_{k=2}^n\frac{{}_0F_1\left(\alpha+1;\mu \lambda_k\right)}{(\lambda_k-x)\displaystyle\prod_{\substack{i=2\\i\neq k}}^n\left(\lambda_k-\lambda_i\right)}\nonumber\\
& \qquad \qquad \left.\times \prod_{i=2}^n\lambda_i^{\alpha}e^{-\lambda_i}(x-\lambda_i)^2
\Delta_{n-1}^2(\boldsymbol{\lambda}){\rm d}\lambda_2\cdots{\rm d}\lambda_n\right) {\rm d}x.
\end{align}
The remainder of this proof is focused on evaluating the above two multiple integrals. Since the two integrals do not share a common structure, in what follows, we will evaluate them separately. 

Let us begin with (\ref{Pdef}). Clearly, the inner multiple integral is symmetric in $\lambda_2,\cdots,\lambda_n$. Thus, we can remove the ordered region of integration to yield
\begin{align*}
\mathfrak{P}(s)=\frac{\mathcal{K}_{n,\alpha}}{(n-1)!}
\frac{e^{-\mu-s}}{\mu^{n-1}}&\int_0^\infty
e^{-x}x^\alpha{}_0F_1\left(\alpha+1;\mu x\right)\nonumber\\
& \qquad \times \left(\int_{[x,\infty)^{n-1}}
\prod_{i=2}^n e^{-\left(1+\frac{s}{x}\right)\lambda_i}\lambda_i^\alpha (x-\lambda_i)
\Delta_{n-1}^2(\boldsymbol{\lambda}){\rm d}\lambda_2\cdots{\rm d}\lambda_n\right)\; {\rm d}x.
\end{align*}
Now we apply the change of variables $y_{i-1}=\frac{(x+s)}{x}(\lambda_i-x),\; i=2,3,\cdots,n$, to the inner $(n-1)$ fold integral with some algebraic manipulation to obtain
\begin{align*}
\mathfrak{P}(s)=(-1)^{(n-1)(1+\alpha)}\frac{\mathcal{K}_{n,\alpha}}{(n-1)!}
 \frac{e^{-\mu-sn}}{\mu^{n-1}} \int_0^\infty
e^{-nx}x^{n(n-1+\alpha)}&\frac{{}_0F_1\left(\alpha+1;\mu x\right)}{(x+s)^{(n+\alpha)(n-1)}}\nonumber\\
& \qquad \times R_{n-1}(-(s+x),\alpha) {\rm d}x
\end{align*}
where
\begin{align*}
R_n(a,\alpha)=\int_{[0,\infty)^n}\prod_{i=1}^n e^{-y_j}y_j(a-y_j)^\alpha \Delta^2_n(\mathbf{y}){\rm d}y_1{\rm d}y_2\cdots{\rm d}y_n.
\end{align*}
Following \cite[Section 22.2.2]{Mehta}, we can solve the above integral to yield\footnote{Specific steps pertaining to this evaluation are not given here as the detailed steps of solving an analogous integral have been given in \cite{PrathaSIAM}.}
\begin{align}
R_n(a,\alpha)=(-1)^{n\alpha}\prod_{j=0}^{n-1}(j+1)!(j+1)!\prod_{j=0}^{\alpha-1}\frac{(n+j)!}{j!}
\det\left[L^{(j)}_{n+i-j}(a)\right]_{i,j=1,2,\cdots,\alpha}.
\end{align}
Therefore, we obtain
\begin{align}
\label{Ppartans}
\mathfrak{P}(s)&=(-1)^{(n-1)}\frac{(n-1)!}{\alpha!}
 \frac{e^{-\mu-sn}}{\mu^{n-1}} \nonumber\\
 &\qquad \times \int_0^\infty
e^{-nx}x^{n(n-1+\alpha)}\frac{{}_0F_1\left(\alpha+1;\mu x\right)}{(x+s)^{(n+\alpha)(n-1)}}
 \det\left[L^{(j)}_{n+i-j-1}(-x-s)\right]_{i,j=1,2,\cdots,\alpha}
{\rm d}x.
\end{align}
Although further manipulation in this form is feasible, it is convenient to leave the solution in the current form. Next we focus on solving the multiple integral given (\ref{Sdef}).

By symmetry, we convert the ordered region of integration in (\ref{Sdef})  to an unordered region to yield
\begin{align*}
\mathfrak{S}(s)=\frac{\mathcal{K}_{n,\alpha}}{(n-1)!}
\frac{e^{-\mu-s}}{\mu^{n-1}}
\int_0^\infty e^{-x}x^{\alpha} &\Biggl(\int_{[x,\infty)^{n-1}}
\sum_{k=2}^n\frac{{}_0F_1\left(\alpha+1;\mu \lambda_k\right)}{(\lambda_k-x)\displaystyle\prod_{\substack{i=2\\i\neq k}}^n\left(\lambda_k-\lambda_i\right)}\nonumber\\
& \qquad \quad \left.\times \prod_{i=2}^n\lambda_i^{\alpha}e^{-\lambda_i}(x-\lambda_i)^2
\Delta_{n-1}^2(\boldsymbol{\lambda}){\rm d}\lambda_2\cdots{\rm d}\lambda_n\right)\; {\rm d}x.
\end{align*}
Since each term in the above summation contributes the same amount, we can further simplify the multiple integral giving
\begin{align*}
\mathfrak{S}(s)=\frac{\mathcal{K}_{n,\alpha}}{(n-2)!}
\frac{e^{-\mu-s}}{\mu^{n-1}}
\int_0^\infty e^{-x}x^{\alpha} &\Biggl(\int_{[x,\infty)^{n-1}}
\frac{{}_0F_1\left(\alpha+1;\mu \lambda_2\right)}{(\lambda_2-x)\displaystyle\prod_{i=3}^n\left(\lambda_2-\lambda_i\right)}\nonumber\\
& \times \prod_{i=2}^n\lambda_i^{\alpha}e^{-\lambda_i}(x-\lambda_i)^2
\Delta_{n-1}^2(\boldsymbol{\lambda}){\rm d}\lambda_2\cdots{\rm d}\lambda_n\Biggr)\; {\rm d}x,
\end{align*}
from which we obtain after using the decomposition $\Delta_{n-1}^2(\boldsymbol{\lambda})=\prod_{j=3}^n(\lambda_2-\lambda_j)^2\Delta_{n-2}^2(\boldsymbol{\lambda})$,
\begin{align*}
\mathfrak{S}(s)
&=\frac{\mathcal{K}_{n,\alpha}}{(n-2)!}
\frac{e^{-\mu-s}}{\mu^{n-1}}
\int_0^\infty e^{-x}x^{\alpha}\left\{\int_x^\infty \lambda_2^\alpha (\lambda_2-x) {}_0F_1\left(\alpha+1;\mu \lambda_2\right)e^{-\left(1+\frac{s}{x}\right)\lambda_2}\right.\nonumber\\
& \times\left. \left(
\int_{[x,\infty)^{n-2}}
\prod_{i=3}^n\lambda_i^{\alpha}e^{-\left(1+\frac{s}{x}\right)\lambda_i}(\lambda_2-\lambda_i)(x-\lambda_i)^2
\Delta_{n-2}^2(\boldsymbol{\lambda}){\rm d}\lambda_3\cdots{\rm d}\lambda_n\right) {\rm d}\lambda_2\right\} {\rm d}x.
\end{align*}
Now we apply the variable transformations
\begin{align*}
y &=\lambda_2-x\\
y_{i-2}&=\frac{(x+s)}{x}(\lambda_i-x),\;\; i=3,4,\cdots,n
\end{align*}
in the above multiple integral to yield
\begin{align*}
 \mathfrak{S}(s)
&=(-1)^{n\alpha}\frac{\mathcal{K}_{n,\alpha}}{(n-2)!}
\frac{e^{-\mu-sn}}{\mu^{n-1}}
\int_0^\infty \frac{e^{-xn}x^\alpha}{\left(1+\frac{s}{x}\right)^{(n-2)(n+\alpha+1)}} \left\{\int_0^\infty y (y+x)^\alpha e^{-\left(1+\frac{s}{x}\right)y}
\right.\nonumber\\
& \hspace{3.5cm} \times{}_0F_1\left(\alpha+1;\mu( y+x)\right) T_{n-2}\left(y\left(1+\frac{s}{x}\right),-s-x,\alpha\right) {\rm d}y\Biggr\}
 {\rm d}x
\end{align*}
where
\begin{align}
T_n(a,b,\alpha):=\int_{[0,\infty)^n}\prod_{i=1}^n(a-y_i)(b-y_i)^\alpha e^{-y_i}y_i^2 \Delta_n^2(\mathbf{y}) {\rm d}y_1{\rm d}y_2\cdots{\rm d}y_n.
\end{align}
It is not difficult to observe that $T_n(a,b,\alpha)$ and $Q_n(a,b,\alpha)$ defined in (\ref{Qintdef}) share a common structure up to a certain Laguerre weight. Therefore, we can readily follow similar arguments as shown in the Appendix with the modified monic orthogonal polynomials given by $\mathsf{P}_k(x)=(-1)^k k! L_k^{(2)}(x)$ to arrive at
\begin{align*}
T_n(a,b,\alpha):=\frac{(-1)^{n+\alpha(n+\alpha)}\widetilde{\mathcal{K}}_{n,\alpha}}{(b-a)^\alpha}\det\left[L^{(2)}_{n+i-1}(a)\;\;\; L_{n+i+1-j}^{(j)}(b)\right]_{\substack{i=1,2,\cdots,\alpha+1\\j=2,3,\cdots,\alpha+1}}
\end{align*}
where
\begin{align}
\widetilde{\mathcal{K}}_{n,\alpha}=\frac{\prod_{j=1}^{\alpha+1}(n+j-1)! \prod_{j=0}^{n-1}(j+1)!(j+2)!}{\prod_{j=0}^{\alpha-1}j!}.
\end{align}
This in turn gives
\begin{align*}
\mathfrak{S}(s)
=(-1)^{n}\frac{(n-1)!}{\alpha!}
\frac{e^{-\mu-sn}}{\mu^{n-1}}&
\int_0^\infty \frac{e^{-xn}x^{\alpha}}{\left(1+\frac{s}{x}\right)^{(n-1)(n+\alpha)-2}} \left\{\int_0^\infty y  e^{-\left(1+\frac{s}{x}\right)y}
{}_0F_1\left(\alpha+1;\mu( y+x)\right) \right.\nonumber\\
&  \times\det\left[L^{(2)}_{n+i-3}\left(y\left(1+\frac{s}{x}\right)\right)\;\;\; L_{n+i-1-j}^{(j)}\left(-x-s\right)\right]_{\substack{i=1,2,\cdots,\alpha+1\\j=2,3,\cdots,\alpha+1}} {\rm d}y\Biggr\}
 {\rm d}x
\end{align*}
from which we obtain after the variable transformation $y\left(1+\frac{s}{x}\right)=t$
%\begin{align}
%& \mathfrak{S}(s)\nonumber\\
%&=(-1)^{n}\widetilde{\mathcal{K}}_{n-2,\alpha}\frac{\mathcal{K}_{n,\alpha}}{(n-2)!}
%\frac{e^{-\mu-sn}}{\mu^{n-1}}
%\int_0^\infty \frac{e^{-xn}x^{\alpha}}{\left(1+\frac{s}{x}\right)^{(n-1)(n+\alpha)}} \left\{\int_0^\infty t  e^{-t}
%\right.\nonumber\\
%& \qquad \quad \times{}_0F_1\left(\alpha+1;\mu\left(x+\frac{t}{1+\frac{s}{x}}\right)\right) \det\left[L^{(2)}_{n+i-3}(t)\;\;\; L_{n+i-1-j}^{(j)}\left(-x-s\right)\right]_{\substack{i=1,2,\cdots,\alpha+1\\j=2,3,\cdots,\alpha+1}} {\rm d}t\Biggr\}
% {\rm d}x.
%\end{align}
%Noting that only the first column of the determinant depends through the variable $t$, we can obtain
\begin{align}
\label{demmel_decom}
\mathfrak{S}(s)
=(-1)^{n}\frac{(n-1)!}{\alpha!}
\frac{e^{-\mu-sn}}{\mu^{n-1}}&
\int_0^\infty \frac{e^{-xn}x^{\alpha}}{\left(1+\frac{s}{x}\right)^{(n-1)(n+\alpha)}} 
\nonumber\\
& \times \det\left[\varrho_i(s,x) \;\;\; L_{n+i-1-j}^{(j)}\left(-x-s\right)\right]_{\substack{i=1,2,\cdots,\alpha+1\\j=2,3,\cdots,\alpha+1}} 
 {\rm d}x
\end{align}
where
\begin{align}
\label{aux_def}
\varrho_i(s,x)=\int_0^\infty t  e^{-t}  {}_0F_1\left(\alpha+1;\mu\left(x+\frac{t}{1+\frac{s}{x}}\right)\right)  L^{(2)}_{n+i-3}(t){\rm d}t 
\end{align}
and we have used the fact that only the first column of the determinant depends through the variable $t$. The integral in (\ref{aux_def}) does not seem to have a simple closed form solution. Therefore, to facilitate further analysis, we write the hypergeometric function with its equivalent power series expansion and use Lemma \ref{lag} to arrive at
%\begin{align*}
%& \mathfrak{S}(s)\nonumber\\
%%&=(-1)^{n}\widetilde{\mathcal{K}}_{n-2,\alpha}\frac{\mathcal{K}_{n,\alpha}}{(n-2)!}
%%\frac{e^{-\mu-sn}}{\mu^{n-1}}
%%\int_0^\infty \frac{e^{-xn}x^{\alpha}}{\left(1+\frac{s}{x}\right)^{(n-1)(n+\alpha)}} 
%%\nonumber\\
%%& \qquad \quad \times
%%\det\left[\sum_{k=0}^\infty \frac{\mu ^k}{k! (\alpha+1)_k}\int_0^\infty t  e^{-t}\left(x+\frac{t}{1+\frac{s}{x}}\right)^k  L^{(2)}_{n+i-3}(t){\rm d}t  \;\;\; L_{n+i-1-j}^{(j)}\left(-x-s\right)\right]_{\substack{i=1,2,\cdots,\alpha+1\\j=2,3,\cdots,\alpha+1}} 
%% {\rm d}x\nonumber\\
% &=(-1)^{n}\widetilde{\mathcal{K}}_{n-2,\alpha}\frac{\mathcal{K}_{n,\alpha}}{(n-2)!}
%\frac{e^{-\mu-sn}}{\mu^{n-1}}
%\int_0^\infty \frac{e^{-xn}x^{\alpha}}{\left(1+\frac{s}{x}\right)^{(n-1)(n+\alpha)}} 
%\nonumber\\
%& \quad \times
%\det\left[\sum_{p=0}^\infty \sum_{k=0}^p\frac{\mu ^p x^{p-k}}{k!(p-k)! (\alpha+1)_p}\frac{1}{\left(1+\frac{s}{x}\right)^k}\int_0^\infty t^{k+1}  e^{-t} L^{(2)}_{n+i-3}(t){\rm d}t  \;\;\; L_{n+i-1-j}^{(j)}\left(-x-s\right)\right]_{\substack{i=1,2,\cdots,\alpha+1\\j=2,3,\cdots,\alpha+1}} \hspace{-2mm}
% {\rm d}x.
%\end{align*}
%Now by invoking Lemma \ref{lag}  we solve the inner integral to obtain
 \begin{align*}
 \varrho_i(s,x)&=\frac{1}{(n+i-3)!}\sum_{p=0}^\infty \sum_{k=0}^p\frac{\mu ^p x^{p-k} (k+1)!}{k!(p-k)! (\alpha+1)_p}\frac{(1-k)_{n+i-3}}{\left(1+\frac{s}{x}\right)^k}\nonumber\\
  &=
  \frac{1}{(n+i-3)!}\sum_{k=0}^\infty \sum_{p=0}^\infty\frac{\mu ^{p+k} x^{p} (k+1)!}{k!p! (\alpha+1)_{p+k}}\frac{(1-k)_{n+i-3}}{\left(1+\frac{s}{x}\right)^k}.
%& \mathfrak{S}(s)\nonumber\\
% &=(-1)^{n}\widetilde{\mathcal{K}}_{n-2,\alpha}\frac{\mathcal{K}_{n,\alpha}}{(n-2)!}
%\frac{e^{-\mu-sn}}{\mu^{n-1}}
%\int_0^\infty \frac{e^{-xn}x^{\alpha}}{\left(1+\frac{s}{x}\right)^{(n-1)(n+\alpha)}} 
%\nonumber\\
%& \quad \times
%\det\left[\frac{1}{(n+i-3)!}\sum_{p=0}^\infty \sum_{k=0}^p\frac{\mu ^p x^{p-k} (k+1)!}{k!(p-k)! (\alpha+1)_p}\frac{(1-k)_{n+i-3}}{\left(1+\frac{s}{x}\right)^k}  \;\;\; L_{n+i-1-j}^{(j)}\left(-x-s\right)\right]_{\substack{i=1,2,\cdots,\alpha+1\\j=2,3,\cdots,\alpha+1}} 
% {\rm d}x
\end{align*}
%from which we obtain after rearranging the summations
%\begin{align}
%& \mathfrak{S}(s)\nonumber\\
% &=(-1)^{n}\widetilde{\mathcal{K}}_{n-2,\alpha}\frac{\mathcal{K}_{n,\alpha}}{(n-2)!}
%\frac{e^{-\mu-sn}}{\mu^{n-1}}
%\int_0^\infty \frac{e^{-xn}x^{\alpha}}{\left(1+\frac{s}{x}\right)^{(n-1)(n+\alpha)}} 
%\nonumber\\
%& \quad \times
%\det\left[\frac{1}{(n+i-3)!}\sum_{l=0}^\infty \sum_{k=l}^\infty\frac{\mu ^k x^{k-l} (l+1)!}{l!(k-l)! (\alpha+1)_k}\frac{(1-l)_{n+i-3}}{\left(1+\frac{s}{x}\right)^l}  \;\;\; L_{n+i-1-j}^{(j)}\left(-x-s\right)\right]_{\substack{i=1,2,\cdots,\alpha+1\\j=2,3,\cdots,\alpha+1}} 
% {\rm d}x.
%\end{align}
%A simple change of the summation index, $k-l=p$, in the inner summation then yields
%\begin{align*}
%& \mathfrak{S}(s)\nonumber\\
% &=(-1)^{n}\widetilde{\mathcal{K}}_{n-2,\alpha}\frac{\mathcal{K}_{n,\alpha}}{(n-2)!}
%\frac{e^{-\mu-sn}}{\mu^{n-1}}
%\int_0^\infty \frac{e^{-xn}x^{\alpha}}{\left(1+\frac{s}{x}\right)^{(n-1)(n+\alpha)}} 
%\nonumber\\
%& \quad \times
%\det\left[\frac{1}{(n+i-3)!}\sum_{k=0}^\infty \sum_{p=0}^\infty\frac{\mu ^{p+k} x^{p} (k+1)!}{k!p! (\alpha+1)_{p+k}}\frac{(1-k)_{n+i-3}}{\left(1+\frac{s}{x}\right)^k}  \;\;\; L_{n+i-1-j}^{(j)}\left(-x-s\right)\right]_{\substack{i=1,2,\cdots,\alpha+1\\j=2,3,\cdots,\alpha+1}} 
% {\rm d}x.
%\end{align*}
The behavior of the Pochhammer symbol $(1-k)_{n+i-3}$ with respect to $l$ deserves a special attention at this juncture. As such, we can observe that 
\begin{equation*}
(1-k)_{n+i-3}=\left\{\begin{array}{cc}
(n+i-3)! & \text{for $k=0$}\\
0 & \text{for $k=1$}\\
(1-k)_{n+i-3} & \text{for $k\geq 2$},
\end{array}\right.
\end{equation*}
which enables us to decompose the terms corresponding to the summation index $k$ into two parts. As a result, after some algebra,  we obtain
\begin{align}
\label{sigma_decom}
\varrho_i(s,x)&=
{}_0F_1\left(\alpha+1;\mu x\right) +\frac{\sigma_i(s,x)}{(n+i-3)!}.
%
%
%& \mathfrak{S}(s)\nonumber\\
% &=(-1)^{n}\widetilde{\mathcal{K}}_{n-2,\alpha}\frac{\mathcal{K}_{n,\alpha}}{(n-2)!}
%\frac{e^{-\mu-sn}}{\mu^{n-1}}
%\int_0^\infty \frac{e^{-xn}x^{\alpha}}{\left(1+\frac{s}{x}\right)^{(n-1)(n+\alpha)}} 
%\nonumber\\
%& \quad \times
%\det\left[ {}_0F_1\left(\alpha+1;\mu x\right) +\frac{1}{(n+i-3)!}\sum_{k=2}^\infty \sum_{p=0}^\infty\frac{\mu ^{p+k} x^{p} (k+1)!}{k!p! (\alpha+1)_{p+k}}\frac{(1-l)_{n+i-3}}{\left(1+\frac{s}{x}\right)^k}\right.\nonumber\\
%&\hspace{6cm} \;\;\; L_{n+i-1-j}^{(j)}\left(-x-s\right)\Biggr]_{\substack{i=1,2,\cdots,\alpha+1\\j=2,3,\cdots,\alpha+1}} 
% {\rm d}x.
\end{align}
where
\begin{align}
\label{sigma_def}
\sigma_i(s,x)=\sum_{k=0}^\infty \sum_{p=0}^\infty\frac{\mu ^{p+k+2} x^{p} (k+3)!}{(k+2)!p! (\alpha+1)_{p+k+2}}\frac{(-1-k)_{n+i-3}}{\left(1+\frac{s}{x}\right)^{k+2}}.
\end{align}
Now we substitute (\ref{sigma_decom}) into (\ref{demmel_decom}) and further simplify the resultant determinant using the multilinear property to obtain
\begin{align}
\label{multi}
 \mathfrak{S}(s)
 &=(-1)^{n}\frac{(n-1)!}{\alpha!}
\frac{e^{-\mu-sn}}{\mu^{n-1}}
\int_0^\infty e^{-xn}x^{n(n+\alpha-1)}\frac{{}_0F_1\left(\alpha+1;\mu x\right)}{\left(x+s\right)^{(n-1)(n+\alpha)}} \nonumber\\
& \hspace{6.5cm}\times 
\det\left[ 1\;\;\;L_{n+i-1-j}^{(j)}\left(-x-s\right)\right]_{\substack{i=1,2,\cdots,\alpha+1\\j=2,3,\cdots,\alpha+1}} {\rm d}x\nonumber\\
&\qquad 
+(-1)^{n}\frac{(n-1)!}{\alpha!}
\frac{e^{-\mu-sn}}{\mu^{n-1}}\int_0^\infty
\frac{e^{-xn}x^{\alpha}}{\left(1+\frac{s}{x}\right)^{(n-1)(n+\alpha)}}\nonumber\\
& \hspace{5.5cm}\times 
\det\left[\frac{\sigma_i(s,x)}{(n+i-3)!}\;\;\; L_{n+i-1-j}^{(j)}\left(-x-s\right)\right]_{\substack{i=1,2,\cdots,\alpha+1\\j=2,3,\cdots,\alpha+1}} 
 {\rm d}x.
\end{align}
Let us now focus on further simplification of the determinant in the first integral. To this end, we apply the row operation, $i\text{th row}\to i\text{th row}+(-1)(i-1)\text{th row}$ on each row for $i=2,3,\cdots,\alpha+1$ and expand the resultant determinant using its first column to obtain
\begin{align}
\det\left[ 1\;\;\;L_{n+i-1-j}^{(j)}\left(-x-s\right)\right]_{\substack{i=1,2,\cdots,\alpha+1\\j=2,3,\cdots,\alpha+1}}&
%\det\left[ L_{n+i-1-j}^{(j)}\left(-x-s\right)-L_{n+i-2-j}^{(j)}\left(-x-s\right)\right]_{i,j=2,3,\cdots,\alpha+1}\nonumber\\
%& = \det\left[ L_{n+i-1-j}^{(j-1)}\left(-x-s\right)\right]_{i,j=2,3,\cdots,\alpha+1}\nonumber\\
=\det\left[ L_{n+i-1-j}^{(j)}\left(-x-s\right)\right]_{i,j=1,2,\cdots,\alpha}
\end{align}
where we have used the contiguous relation given in (\ref{contg}). Therefore, in view of (\ref{Ppartans}),  we can clearly identify the first term in (\ref{multi}) as $-\mathfrak{P}(s)$. This key observation along with (\ref{mgffact}) gives

\begin{align*}
 \mathfrak{M}_V(s)
=(-1)^{n}\frac{(n-1)!}{\alpha !}
\frac{e^{-\mu-sn}}{\mu^{n-1}} & \int_0^\infty
\frac{e^{-xn}x^{\alpha}}{\left(1+\frac{s}{x}\right)^{(n-1)(n+\alpha)}}\nonumber\\
& \times 
\det\left[\frac{\sigma_i(s,x)}{(n+i-3)!} \;\;\; L_{n+i-1-j}^{(j)}\left(-x-s\right)\right]_{\substack{i=1,2,\cdots,\alpha+1\\j=2,3,\cdots,\alpha+1}} 
 {\rm d}x.
\end{align*}
The remaining task at hand is to further simplify $\sigma_i(s,x)$ given in (\ref{sigma_def}).  To this end,
following (\ref{poch}),  we find that $(-1-k))_{n+i-3}$ is non-zero for $k\geq n+i-4$. Therefore, we shift the index $k$ with some algebraic manipulation to obtain the m.g.f. of $V$ as
\begin{align}
\mathfrak{M}_V(s)&=(n-1)!
e^{-\mu-sn}\int_0^\infty
\frac{e^{-xn}x^{n(n+\alpha)-1}}{\left(x+s\right)^{(n-1)(n+\alpha+1)}}\nonumber\\
& \qquad \qquad \times \det\left[
\left(-\frac{\mu x}{x+s}\right)^{i-1}\vartheta_i(x\mu,x+s)
\;\;\; L_{n+i-1-j}^{(j)}\left(-x-s\right)\right]_{\substack{i=1,2,\cdots,\alpha+1\\j=2,3,\cdots,\alpha+1}} 
 \hspace{-1cm}{\rm d}x
\end{align}
where
\begin{align*}
\vartheta_i(w,z)=\frac{(n+i-1)}{(n+\alpha+i-2)!}
\sum_{k=0}^\infty 
\frac{(n+i)_k(n+i-2)_k\; {}_0F_1\left(\alpha+n+i+k-1;w\right) w^k}{k!(\alpha+n+i-1)_k(n+i-1)_k z^k}.
\end{align*}
Finally, we take the inverse Laplace transform of the above to yield the p.d.f. of $V$ which concludes the proof.

Although further simplification of (\ref{pdfexact}) seems intractable for general matrix dimensions $m$ and $n$, we can obtain a relatively simple expression in the important case of square matrices (i.e., $m=n$), which is given in the following corollary.
\begin{cor}
For $\alpha=0$, (\ref{pdfexact}) becomes
\begin{align}
f^{(0)}_V(v) &=n(n^2-1) e^{-\mu}
(v-n)^{n^2-2}v^{-n^2}\nonumber\\
& \times 
\sum_{k=0}^\infty
\frac{(n^2)_k}{(n)_k k!} \left(\frac{\mu}{v}\right)^k {}_3F_3\left(n+1,n-1,n^2+k;n,n+k,n^2-1;\mu\left(1-\frac{n}{v}\right)\right) H(v-n)
\end{align}
where $H(z)$ denotes the unit step function and ${}_3F_3(a_1,a_2,a_3;c_1,c_2,c_3;z)$ is the generalized hypergeometric function \cite{Erdelyi}.
\end{cor} 

It is also worth pointing out that, for $\mu=0$ (i.e., when the matrix $\mathbf{W}$ is a central Wishart matrix), (\ref{pdfexact}) simplifies to
\begin{align}
f^{(\alpha)}_V(v)=\frac{n!(n^2+n\alpha-1)!}{(n+\alpha-1)! v^{n(n+\alpha)}}
\mathcal{L}^{-1}\left\{\frac{e^{-ns}}{s^{(n-1)(n+\alpha+1)}}\det\left[L^{(j+1)}_{n+i-j-1}(-s)\right]_{i,j=1,2,\cdots,\alpha}\right\}
\end{align}
 which coincides with the corresponding result given in \cite[Corollary 3.2]{PrathaSIAM}.
 
% Figure \ref{Fig2} compares the simulation and analytical p.d.f. results corresponding to the case of $m=n=5$ with $\mu=2$.
%\begin{figure}
% \centering
% \vspace*{1.0cm}
%            \includegraphics[width=.7\textwidth]{5by5mu2.eps}
% \caption{Comparison of simulated data points and the analytical p.d.f.\ $f^{(0)}_V(v)$ with $m=n=5$ and $\mu=2$.}
% \vspace*{0.3cm}
% \label{Fig2}
%\end{figure}

\section{The Average of the Reciprocal of a Certain Characteristic Polynomial}
Here we consider the problem of determining the average of the reciprocal of a certain characteristic polynomial with respect to a complex non-central Wishart density with a rank one mean. It is noteworthy that this particular problem corresponding to complex central Wishart matrices has been solved in \cite{Mehta, Fodorov}. A general framework to derive such averages based on duality relations has been proposed in \cite{Patric}. However, the duality relation given in \cite[Proposition 8]{Patric} does not seem to apply here, since the stringent technical requirements for the validity of that formula are not satisfied by the parameters in our model of interest. Moreover, this particular case has not been considered in a recent detailed analysis on the averages of characteristic polynomials for Gaussian and Chiral Gaussian matrices with an external source \cite{Peter}. Therefore, in what follows, we derive the average of one of the basic forms of the reciprocal of the characteristic polynomial. The most general form, however, is not investigated here.

Let us consider the following average
\begin{equation}
\label{average}
{\rm{E}}_{\mathbf{W}}\left(\frac{1}{\det[z\mathbf{I}_n+\mathbf{W}]}\right)={\rm{E}}_{\boldsymbol{\Lambda}}\left(\prod_{j=1}^n\frac{1}{z+\lambda_j}\right),\; |\arg{z}|<\pi,
\end{equation}
the value of which is given in the following theorem.

\begin{thm}
Let $\mathbf{W}\sim\mathcal{W}_n\left(m,\mathbf{I}_n,\mathbf{M}^\dagger\mathbf{M}\right)$, where $\mathbf{M}$ is rank-$1$ and $\mathrm{tr}(\mathbf{M}^\dagger\mathbf{M})=\mu$. Then (\ref{average}) is given by
\begin{align}
\label{reciprocal}
{\rm{E}}_{\mathbf{W}}\left(\frac{1}{\det[z\mathbf{I}_n+\mathbf{W}]}\right)=
z^\alpha\sum_{k=0}^\infty
(-\mu)^k \Psi(k+n+\alpha;\alpha+1;z),\;\; |\arg{z}|<\pi
\end{align}
where $\Psi(a;c;z)$ is the confluent hypergeometric function of the second kind.
\end{thm}
{\bf{Proof:}} Due to symmetry, we have
\begin{equation*}
{\rm{E}}_{\boldsymbol{\Lambda}}\left(\prod_{j=1}^n\frac{1}{z+\lambda_j}\right)=\frac{1}{n!}\int_{[0,\infty)^n}
\frac{f_{\boldsymbol{\Lambda}}(\lambda_1,\lambda_2,\cdots,\lambda_n)}{\prod_{j=1}^n (z+\lambda_j)} {\rm d}\lambda_1{\rm d}\lambda_2\cdots{\rm d}\lambda_n.
\end{equation*}
Now it is convenient to apply partial faction decomposition to yield
\begin{align*}
{\rm{E}}_{\boldsymbol{\Lambda}}\left(\prod_{j=1}^n\frac{1}{z+\lambda_j}\right)=\frac{1}{n!}\sum_{j=1}^n \int_{[0,\infty)^n}
\frac{1}{\displaystyle \prod_{\substack{i=1\\
i\neq j}}^n(\lambda_i-\lambda_j)}
\frac{f_{\boldsymbol{\Lambda}}(\lambda_1,\lambda_2,\cdots,\lambda_n)}{ (z+\lambda_j)} {\rm d}\lambda_1{\rm d}\lambda_2\cdots{\rm d}\lambda_n.
\end{align*}
Since each integral in the above summation contributes the same amount, we can simplify it as
\begin{align*}
{\rm{E}}_{\boldsymbol{\Lambda}}\left(\prod_{j=1}^n\frac{1}{z+\lambda_j}\right)=\frac{1}{(n-1)!} \int_{[0,\infty)^n}
\frac{1}{\prod_{j=2
}^n(\lambda_j-\lambda_1)}
\frac{f_{\boldsymbol{\Lambda}}(\lambda_1,\lambda_2,\cdots,\lambda_n)}{ (z+\lambda_1)} {\rm d}\lambda_1{\rm d}\lambda_2\cdots{\rm d}\lambda_n.
\end{align*}
To facilitates further analysis, we use (\ref{newden}) with some rearrangements to write
\begin{align}
\label{chaeq}
{\rm{E}}_{\boldsymbol{\Lambda}}\left(\prod_{j=1}^n\frac{1}{z+\lambda_j}\right)=\Omega_1(z)+\Omega_2(z) 
\end{align}
where
\begin{align}
\label{omega1}
\Omega_1(z)&=
%(-1)^{n-1}\frac{\mathcal{K}_{n,\alpha}}{(n-1)!}\frac{e^{-\mu}}{\mu^{n-1}}
%\int_0^\infty
%\frac{{}_0F_1\left(\alpha+1,\mu \lambda_1\right)}{z+\lambda_1} \lambda_1 ^{\alpha} e^{-\lambda_1}\nonumber\\
%& \hspace{3.5cm}\times
%\int_{[0,\infty)^{n-1}}
%\prod_{j=2}^n\frac{\lambda_j^\alpha e^{-\lambda_j}}{(\lambda_j-\lambda_1)^2} \Delta^2_n(\boldsymbol{\lambda}) {\rm d}\lambda_2 {\rm d}\lambda_2\cdots{\rm d}\lambda_n \;  {\rm d}\lambda_1\nonumber\\
\frac{(-1)^{n-1}}{\alpha !}\frac{e^{-\mu}}{\mu^{n-1}}
\int_0^\infty
\frac{{}_0F_1\left(\alpha+1,\mu \lambda_1\right)}{z+\lambda_1} \lambda_1 ^{\alpha} e^{-\lambda_1}{\rm d}\lambda_1
\end{align}
and
\begin{align}
\label{omega2}
\Omega_2(z)& =\frac{\mathcal{K}_{n,\alpha}}{(n-1)!}\frac{e^{-\mu}}{\mu^{n-1}}
\int_0^\infty \frac{\lambda_1^\alpha e^{-\lambda_1}}{z+\lambda_1}\left(
\sum_{k=2}^n
\int_{[0,\infty)^{n-1}}
\frac{{}_0F_1\left(\alpha+1,\mu \lambda_k\right)}{\prod_{\substack{j=1\\ j\neq k}}^n (\lambda_j-\lambda_k)}\right.\nonumber\\
& \hspace{5.5cm}\left.\times \prod_{j=2}^n\frac{\lambda_j^\alpha e^{-\lambda_j}}{ (\lambda_j-\lambda_1)}
\Delta_n^2(\boldsymbol{\lambda}){\rm d}\lambda_2 {\rm d}\lambda_3\cdots{\rm d}\lambda_n \right) {\rm d}\lambda_1.
\end{align}
Since further simplification of (\ref{omega1}) seems an arduous task, we leave it in its current form and focus on (\ref{omega2}). Noting that each term inside the summation contributes the same amount due to symmetry in $\lambda_2,\lambda_3,\cdots,\lambda_n$, we can further simplify (\ref{omega2}) to yield
\begin{align*}
\Omega_2(z)& =\frac{(-1)^n\mathcal{K}_{n,\alpha}}{(n-2)!}\frac{e^{-\mu}}{\mu^{n-1}}
\int_0^\infty \frac{\lambda_1^\alpha e^{-\lambda_1}}{z+\lambda_1}
\int_0^\infty
{}_0F_1\left(\alpha+1,\mu \lambda_2\right) \lambda_2^\alpha e^{-\lambda_2}\left(
\int_{[0,\infty)^{n-1}}
\frac{1}{\prod_{\substack{j=1\\ j\neq 2}}^n (\lambda_j-\lambda_k)}\right.\nonumber\\
& \hspace{7cm}\left.\times \prod_{j=2}^n\frac{\lambda_j^\alpha e^{-\lambda_j}}{ (\lambda_1-\lambda_j)}
\Delta_n^2(\boldsymbol{\lambda}) {\rm d}\lambda_3\cdots{\rm d}\lambda_n \right){\rm d}\lambda_2  {\rm d}\lambda_1.
\end{align*}
We now use the decomposition $\Delta_n^2(\boldsymbol{\lambda})=\prod_{j=2}^n (\lambda_1-\lambda_j)^2\prod_{j=3}^n(\lambda_2-\lambda_j)^2\Delta^2_{n-2}(\boldsymbol{\lambda})$ followed by the variable transformation $y_j=\lambda_{j-2},\;j=3,4,\cdots,n$, to obtain
\begin{align*}
\Omega_2(z)= \frac{(-1)^n\mathcal{K}_{n,\alpha}}{(n-2)!}\frac{e^{-\mu}}{\mu^{n-1}}
\int_0^\infty \frac{\lambda_1^\alpha e^{-\lambda_1}}{z+\lambda_1}\left(
\int_0^\infty
{}_0F_1\left(\alpha+1,\mu \lambda_2\right) \lambda_2^\alpha e^{-\lambda_2}
U_{n-2}(\lambda_1,\lambda_2,\alpha)
{\rm d}\lambda_2 \right) {\rm d}\lambda_1
\end{align*}
where
\begin{align}
U_n(r_1,r_2,\alpha):=\int_{[0,\infty)^{n}} \prod_{j=1}^n\prod_{i=1}^2(r_i-y_j)y^\alpha_j e^{-y_j} \Delta^2_n({\bf{y}}) {\rm d}y_1 {\rm d}y_2\cdots {\rm d}y_n.
\end{align}
The above integral can be solved using \cite[Eqs. 22.4.2, 22.4.11]{Mehta} and the Appendix with the choice of $\mathsf{P}_k(x)=(-1)^kk!L_k^{(\alpha)}(x)$ to yield
\begin{align}
U_n(r_1,r_2,\alpha)=(-1)n!(n+1)!\prod_{j=0}^{n-1}(j+1)!(j+\alpha)!\frac{\det\left[L^{(\alpha)}_{n+i-1}(r_j)\right]_{i,j=1,2}}{(r_2-r_1)}.
\end{align}
This in turn gives
\begin{align}
\label{omega2int}
\Omega_2(z)& = (-1)^{n+1}\frac{(n-1)!}{\alpha ! (n+\alpha-2)!}\frac{e^{-\mu}}{\mu^{n-1}}\nonumber\\
& \times
\int_0^\infty \frac{\lambda_1^\alpha e^{-\lambda_1}}{z+\lambda_1}\left(
\int_0^\infty
{}_0F_1\left(\alpha+1,\mu \lambda_2\right) 
\frac{\det\left[L^{(\alpha)}_{n+i-1}(\lambda_j)\right]_{i,j=1,2}}{(\lambda_2-\lambda_1)}\lambda_2^\alpha e^{-\lambda_2}
{\rm d}\lambda_2 \right) {\rm d}\lambda_1.
\end{align}
Further manipulation of the above integral in its current form is highly undesirable due to the term $\lambda_2-\lambda_1$  in the denominator. To circumvent this difficulty, we employ the following form of the Christoffel-Darboux formula
\begin{align*}
\frac{\det\left[L^{(\alpha)}_{n+i-1}(\lambda_j)\right]_{i,j=1,2}}{(\lambda_2-\lambda_1)}& =
\frac{L^{(\alpha)}_{n-1}(\lambda_2)L^{(\alpha)}_{n-2}(\lambda_1)-L^{(\alpha)}_{n-1}(\lambda_1)L^{(\alpha)}_{n-2}(\lambda_2)}{\lambda_2-\lambda_1}\nonumber\\
&= (-1)\frac{(n+\alpha-2)!}{(n-1)!}\sum_{j=0}^{n-2}\frac{j!}{(j+\alpha)!}L^{(\alpha)}_{j}(\lambda_1)L^{(\alpha)}_{j}(\lambda_2)
\end{align*}
in (\ref{omega2int}) to obtain 
\begin{align*}
\Omega_2(z) = \frac{(-1)^n}{\alpha !}\frac{e^{-\mu}}{\mu^{n-1}}
\sum_{j=0}^{n-2}\frac{j!}{(j+\alpha)!}& 
\int_0^\infty \frac{L^{(\alpha)}_{j}(\lambda_1)}{z+\lambda_1}\lambda_1^\alpha e^{-\lambda_1}{\rm d}\lambda_1\\
&\qquad \quad \times \int_0^\infty
{}_0F_1\left(\alpha+1,\mu \lambda_2\right) L^{(\alpha)}_{j}(\lambda_2)\lambda_2^\alpha e^{-\lambda_2} {\rm d}\lambda_2.
\end{align*}
The second integral can be solved using Lemma \ref{lag} to obtain
\begin{align}
\int_0^\infty
{}_0F_1\left(\alpha+1,\mu \lambda_2\right) L^{(\alpha)}_{j}(\lambda_2)\lambda_2^\alpha e^{-\lambda_2} {\rm d}\lambda_2=
\frac{\alpha!}{j!}(-\mu)^je^{\mu}
\end{align}
which in turn gives
\begin{align*}
\Omega_2(z)= \frac{(-1)^{n}}{\alpha !}\frac{e^{-\mu}}{\mu^{n-1}}
\sum_{j=0}^{n-2}\frac{(-\mu)^j\alpha !}{(j+\alpha)!}\;e^{\mu}
\int_0^\infty \frac{L^{(\alpha)}_{j}(\lambda_1)}{z+\lambda_1}\lambda_1^\alpha e^{-\lambda_1}{\rm d}\lambda_1. 
\end{align*}
In order to further simplify the above integral, we rearrange the summation with respect to index $j$ giving
\begin{align}
\label{omega2int1}
\Omega_2(z) = \frac{(-1)^{n}}{\alpha !}\frac{e^{-\mu}}{\mu^{n-1}} & \left(
\sum_{j=0}^{\infty}\frac{(-\mu)^j\alpha !}{(j+\alpha)!}\;e^{\mu}
\int_0^\infty \frac{L^{(\alpha)}_{j}(\lambda_1)}{z+\lambda_1}\lambda_1^\alpha e^{-\lambda_1}{\rm d}\lambda_1\right.\nonumber\\
& \hspace{2.6cm}-\left.
\sum_{j=n-1}^{\infty}\frac{(-\mu)^j\alpha !}{(j+\alpha)!}\;e^{\mu}
\int_0^\infty \frac{L^{(\alpha)}_{j}(\lambda_1)}{z+\lambda_1}\lambda_1^\alpha e^{-\lambda_1}{\rm d}\lambda_1\right).
\end{align}
Let us now focus on the first infinite summation. As such, using (\ref{lagdef}) with some algebraic manipulation we get
\begin{align}
\sum_{j=0}^{\infty}\frac{(-\mu)^j\alpha !}{(j+\alpha)!}L^{(\alpha)}_{j}(\lambda_1)& ={}_0F_1(\alpha+1;\mu \lambda_1) e^{-\mu}.
\end{align}
%from which we obtain after shifting the initial value of the summation index $j$ to zero and re-summing the resultant series
%\begin{align}
%\sum_{j=0}^{\infty}\frac{(-\mu)^j\alpha !}{(j+\alpha)!}L^{(\alpha)}_{j}(\lambda_1) =
%\end{align} 
Therefore, (\ref{omega2int1}) simplifies to
\begin{align*}
\Omega_2(z)= \frac{(-1)^{n}}{\alpha !}\frac{e^{-\mu}}{\mu^{n-1}}&\left(
\int_0^\infty \frac{{}_0F_1(\alpha+1;\mu \lambda_1)}{z+\lambda_1}\lambda_1^\alpha e^{-\lambda_1}{\rm d}\lambda_1\right.\\
&\hspace{3cm}-\left.
\sum_{j=n-1}^{\infty}\frac{(-\mu)^j\alpha !}{(j+\alpha)!}\;e^{\mu}
\int_0^\infty \frac{L^{(\alpha)}_{j}(\lambda_1)}{z+\lambda_1}\lambda_1^\alpha e^{-\lambda_1}{\rm d}\lambda_1\right).
\end{align*}
from which, in view of  (\ref{omega1}), we obtain
\begin{align*}
\Omega_2(z)& =-\Omega_1(z)+ (-1)^{n+1}\frac{1}{\mu^{n-1}}
\sum_{j=n-1}^{\infty}\frac{(-\mu)^j}{(j+\alpha)!}
\int_0^\infty \frac{L^{(\alpha)}_{j}(\lambda_1)}{z+\lambda_1}\lambda_1^\alpha e^{-\lambda_1}{\rm d}\lambda_1.
\end{align*}
Finally, we shift the initial value of the summation index to zero and use \cite[Eq. 6.15.2.16]{Erdelyi} with (\ref{chaeq}) to yield (\ref{reciprocal}) which concludes the proof.

\section*{Acknowledgment}
The author would like to thank Yang Chen and Matthew McKay for insightful discussions. This work was supported by a National Science Foundation grant.
\section*{Appendix}
Following \cite[Eqs. 22.4.2, 22.4.11]{Mehta}, we begin with the integral
\begin{align}
\label{qbeg}
\int_{[0,\infty)^n}\prod_{j=1}^ne^{-y_j}\prod_{i=1}^{\alpha+1}(r_i-y_j)&\Delta_n^2(\mathbf{y}) {\rm d}y_1 {\rm d}y_2\cdots{\rm d}y_n = \prod_{i=0}^{n-1}(i+1)!i!\;\frac{\det\left[\mathsf{P}_{n+i-1}(r_j)\right]_{i,j=1,2,\cdots,\alpha+1}}{\Delta_{\alpha+1}(\mathbf{r})},
\end{align}
where $\mathsf{P}_{k}(x)$'s are monic polynomials orthogonal with respect to $e^{-x}$, over $0\leq x<\infty$. As such, we choose $\mathsf{P}_k(x)=(-1)^kk!L_k^{(0)}(x)$, which upon substituting into the above equation gives
\begin{align}
& \int_{[0,\infty)^n}\prod_{j=1}^ne^{-y_j}\prod_{i=1}^{\alpha+1}(r_i-y_j)\Delta_n^2(\mathbf{y}) {\rm d}y_1 {\rm d}y_2\cdots{\rm d}y_n \nonumber\\
&=(-1)^{(n-1)(\alpha+1)}\prod_{i=0}^{n-1}(i+1)!i!\prod_{i=1}^{\alpha+1}(-1)^i(n+i-1)!\; \frac{\det\left[L^{(0)}_{n+i-1}(r_j)\right]_{i,j=1,2,\cdots,\alpha+1}}{\Delta_{\alpha+1}(\mathbf{r})}.
\end{align}
In general, the $r_i$'s in the above formula are distinct parameters. However, for our purpose, we have to choose them in such a manner that the left side of (\ref{qbeg}) becomes $Q_n(a,b,\alpha)$. To this end, we select $r_i$'s such that
\begin{equation*}
r_i=\left\{\begin{array}{ll}
a & \text{if $i=1$}\\
b & \text{if $i=2,3,\cdots,\alpha+1$.}
\end{array}\right.
\end{equation*}
This direct substitution in turn gives a $\frac{0}{0}$ indeterminate form for the right side of (\ref{qbeg}). To circumvent this problem, instead of direct substitution, we evaluate the following limit
\begin{align}
\label{qdef}
Q_n(a,b,\alpha)&=(-1)^{(n-1)(\alpha+1)}\prod_{i=0}^{n-1}(i+1)!i!\prod_{i=1}^{\alpha+1}(-1)^i(n+i-1)!\nonumber\\
& \hspace{2cm} \times \lim_{r_2,r_3,\cdots,r_{\alpha+1}\to b}\frac{\det\left[L^{(0}_{n+i-1}(a)\;\;\;L^{(0)}_{n+i-1}(r_j)\right]_{\substack{i=1,2,\cdots,\alpha+1\\
j=2,3,\cdots,\alpha+1}}}{\det[a^{i-1}\;\;\; r_j^{i-1}]_{\substack{i=1,2,\cdots,\alpha+1\\
j=2,3,\cdots,\alpha+1}}}.
\end{align}
The limit on the right can be evaluated based on an approach given in \cite{Khatri} to yield
\begin{align}
\label{qlim}
&\lim_{r_2,r_3,\cdots,r_{\alpha+1}\to b}\frac{\det\left[L^{(0)}_{n+i-1}(a)\;\;\;L^{(0)}_{n+i-1}(r_j)\right]_{\substack{i=1,2,\cdots,\alpha+1\\
j=2,3,\cdots,\alpha+1}}}{\det[a^{i-1}\;\;\; r_j^{i-1}]_{\substack{i=1,2,\cdots,\alpha+1\\
j=2,3,\cdots,\alpha+1}}}\nonumber\\
&\hspace{5.5cm}=\frac{\det\left[L^{(0)}_{n+i-1}(a)\;\;\;\displaystyle \frac{{\rm d}^{j-2}}{{\rm d}b^{j-2}}L^{(0)}_{n+i-1}(b)\right]_{\substack{i=1,2,\cdots,\alpha+1\\
j=2,3,\cdots,\alpha+1}}}{\det\left[a^{i-1}\;\;\; \displaystyle \frac{{\rm d}^{j-2}}{{\rm d}b^{j-2}}b^{i-1}\right]_{\substack{i=1,2,\cdots,\alpha+1\\
j=2,3,\cdots,\alpha+1}}}.
\end{align}
The denominator of (\ref{qlim}) gives
\begin{equation}
\label{denom}
\det\left[a^{i-1}\;\;\; \displaystyle \frac{{\rm d}^{j-2}}{{\rm d}b^{j-2}}b^{i-1}\right]_{\substack{i=1,2,\cdots,\alpha+1\\
j=2,3,\cdots,\alpha+1}}=\prod_{i=1}^{\alpha-1} i!\;(b-a)^\alpha.
\end{equation}
The numerator can be simplified using (\ref{lagderi}) to yield
\begin{align}
\label{num}
&\det\left[L^{(0)}_{n+i-1}(a)\;\;\;\displaystyle \frac{{\rm d}^{j-2}}{{\rm d}b^{j-2}}L^{(0)}_{n+i-1}(b)\right]_{\substack{i=1,2,\cdots,\alpha+1\\
j=2,3,\cdots,\alpha+1}}\nonumber\\
&\hspace{3.5cm} = 
(-1)^{\frac{1}{2}\alpha(\alpha-1)}\det\left[L^{(0}_{n+i-1}(a)\;\;\;L^{(j-2)}_{n+i+1-j}(b)\right]_{\substack{i=1,2,\cdots,\alpha+1\\
j=2,3,\cdots,\alpha+1}}.
\end{align}
Substituting (\ref{denom}) and (\ref{num}) into (\ref{qlim}) and then the result into (\ref{qdef}) gives (\ref{q1}).


\begin{thebibliography}{1}

\bibitem{Anderson}
T. W. Anderson, \emph{An Introduction to Multivariate Statistical Analysis}, John Wiley \& Sons Inc., 2003.

\bibitem{Askey}
 G. E. Andrews, R. Askey, and R. Roy, \emph{Special Functions}, Cambridge University Press, Cambridge, U. K., 1999.
 
 \bibitem{Forrester1}
K. E. Bassler, P. J. Forrester, and N. E. Frankel, \emph{Eigenvalue separation in some random matrix models}, arXiv: 0810.1554v1, 2008.

\bibitem{Arous}
 G. Ben Arous and S. P{\'e}ch{\'e}, \emph{Universality of local eigenvalue statistics for some sample covariance matrices}, Commun. Pure and Appl. Math., 58 (2005), pp. 1316--1357.


\bibitem{Bleher2}
P. M. Bleher and A. Kuijlaars, \emph{Random matrices with external source and multiple orthogonal polynomials}, Int. Math. Res. Notice, pp. 109--129, 2004.

\bibitem{Bleher1}
P. M. Bleher and A. Kuijlaars, \emph{Large $n$ limit of Gaussian random matrices with external source, part I}, Comm. Math. Phys. 252 (2004), pp. 43--76.

\bibitem{Bleher3}
P. M. Bleher and A. Kuijlaars, \emph{Integral representations for multiple Hermitian and multiple Laguerre polynomials}, Annales de l`Institut Fourier 55 (2005), pp. 2001--2014.

\bibitem{Boro}
A. Borodin and E. Strahov, \emph{Averages of  characteristic polynomials in random matrix theory }, Commun. Pure Appl. Math., 59 (2006), pp. 161--253.

\bibitem{Brezin1}
E. Br\'ezin and S. Hikami, \emph{Correlations of nearby levels induced by a random potential}, Nucl. Phys. B 479 (1996), pp. 697--706.

\bibitem{Brezin2}
E. Br\'ezin and S. Hikami, \emph{Spectral form factor in random matrix theory}, Phys. Rev. E 55 (1997), pp. 4067--4083.

\bibitem{Brezin4}
E. Br\'ezin and S. Hikami, \emph{Level spacing of random matrices in an external source}, Phys. Rev. E 58 (1998), pp. 7176--7185. 

\bibitem{Burg}
 P. B\"urgisser, F. Cucker, and M. Lotz, \emph{Smoothed analysis of complex condition numbers}, J. Math. Pures. Appl., 86 (2006), pp. 293--309.

\bibitem{Profchen}
Y. Chen, D.-Z. Liu, and D.-S. Zhou, \emph{Smallest eigenvalue distribution of the fixed-trace Laguerre beta ensemble}, J. Phys. A 43 (2010), pp. 1--12.

\bibitem{Cons}
G. Constantine,  \emph{Some noncentral distribution problems in multivariate analysis}, Ann. Math. Statist., 34 (1963), pp. 1270-1285.

\bibitem{Cucker}
 F. Cucker, H. Diao, and Y. Wei, \emph{Smoothed analysis of some conditioned numbers}, Numer. Linear Algebra  Appl., 13 (2006), pp. 71--84.
 
 
\bibitem{Demmel}
J. W. Demmel, \emph{The probability that a numerical analysis problem is difficult}, Math. Comp., 50 (1988), pp. 449--480. 


\bibitem{Patric}
P. Desrosiers, \emph{Duality in random matrix ensembles for all $\beta$}, Nucl. Phys. B 817 (2009), pp. 224--251.

\bibitem{Pratha}
 P. Dharmawansa and M. R. McKay, \emph{Extreme eigenvalue distributions of some complex correlated non-central Wishart and gamma-Wishart random matrices}, JMVA, 102 (2011), pp. 847--868.

\bibitem{PrathaSIAM}
 P. Dharmawansa, M. R. McKay, and Y. Chen, \emph{Distributions of Demmel and related condition numbers}, SIAM J. Matrix Anal. Appl., 34 (2013), pp. 257--279.

\bibitem{Erdelyi}
 A. Erd\'elyi, \emph{Higher Transcendental Functions}, vol. 1, McGraw-Hill, 1953.

\bibitem{Forrester}
 P. J. Forrester and T. D. Hughes, \emph{Complex Wishart matrices and conductance in mesoscopic systems: Exact results}, J. Math. Phys. 35 (1994), pp. 6736--6747.

\bibitem{PJF}
P. J. Forrester, \emph{Log Gases and Random Matrices}, Princeton University Press, Princeton, NJ, 2010. 

\bibitem{Peter}
P. J. Forrester, \emph{The averaged characteristic polynomial for the Gaussian and chiral Gaussian ensembles with a source},  arXiv: 1203.5838v1, 2012.

\bibitem{Fodorov}
Y. V. Fyodorov and E. Strahov, \emph{An exact formula for general spectral correlation function of Hermitian matrices}, J. Phys. A: Math. Gen., 36 (2003), pp. 3203--3213.

\bibitem{Grad}
 I. S. Gradshteyn and I. M. Ryzhik, \emph{Tables of Integrals Series and Products}, 5th ed., Academic Press, New York, 1994


\bibitem{Herz}
C. S. Herz, \emph{Bessel functions of matrix argument}, Ann. Math., 61 (1955), pp. 474--523.


\bibitem{Hua}
L.-K. Hua, \emph{Harmonic Analysis of Functions of Several Complex Variables in the Classical Domains}, Moscow (English translation: Amer. Math. Soc., Providence, R. I.), 1959.

\bibitem{James}
 A. T. James, \emph{Distributions of matrixvariates and latentroots derived from normal samples}, Ann. Math. Stat., 35 (1964), pp. 475--501.

\bibitem{Imj}
I. M. Johnstone, \emph{On the distribution of the largest principal component}, Ann. Math. Statist., 29 (2001), pp. 295--327.

\bibitem{Kazakov}
V. Kazakov, \emph{External matrix field problem and new multicriticalities in $(-2)$- dimensional random surfaces}, Nucl. Phys. B 354 (1991), pp. 614--624.

\bibitem{Khatri}
C. G. Khatri, \emph{On the moments of traces of two matrices in three situations for complex multivariate normal populations}, Sankhy$\bar{a}$, 32 (1970), pp. 65--80.

\bibitem{Krish}
P. R. Krishnaiah and F. J. Schuurmann, \emph{On the evaluation some distributions that arise in simultaneous tests for the equality of the latent roots of the covariance matrix}, J. Multivariate Anal., 4 (1974), pp. 265--282.


\bibitem{Matthesis}
M. R. McKay, \emph{Random matrix theory analysis of multiple-antenna communication systems}, Ph.D. Dissertation, University of Sydney, 2006.


\bibitem{Mehta} 
M. L. Mehta, \emph{Random Matrices}, 3rd ed., Academic Press, New York, 2004.

\bibitem{Mo}
M. Y. Mo,\emph{The rank $1$ real Wishart spiked model I. Finite $N$ analysis}, arXiv: 1011.5404v1, 2010.

\bibitem{Robb} 
R. J. Muirhead, \emph{Aspects of Multivariate Statistical Theory}, Wiley, 1982.

\bibitem{Alex}
A. Onatski, M. J. Moreira, and M. Hallin, \emph{Asymptotic power of sphericity tests for high-dimensional data}, Ann. Math. Stat., 34 (2013), pp. 1204--1231.

\bibitem{Sankar}
 A. Sankar, D. A. Spielman, and S.-H. Teng, \emph{Smoothed analysis of the condition numbers and growth factors of matrices}, SIAM J. Matrix Anal. Appl., 28 (206), pp. 446--476.
 
\bibitem{Spielman}
 D. A. Spielman and S.-H. Teng, \emph{Smoothed analysis: Why the simplex algorithm usually takes polynomial time}, J. ACM, 51 (2004), pp. 385--463.

\bibitem{Szego}
G. Szeg\"o, \emph{Orthogonal Polynomials}, 4th ed., AMS Colloquium Publications, 1975.

\bibitem{Takemura}
 A. Takemura, \emph{Zonal Polynomials}, IMS Lecture Notes-Monograph Series, vol. 4, 1984. 


\bibitem{Telatar}
E. Telatar,\emph{Capacity of Multi-antenna Gaussian channels}, Eur. Tran. Tel., 10 (1999), pp. 585--595.


\bibitem{Verdu}
A. M. Tulino and S. Verd\'u, \emph{Random Matrix Theory and Wireless Communications}, Found. Trend. Commun. Inf. Theory, 2004.


\bibitem{Wang}
 D. Wang, \emph{The largest eigenvalue of real symmetric, Hermitian and Hermitian self-dual random matrix models with rank one external source, part I}, arXiv:1012.4144v3.


 \bibitem{Mario}
 M. Wschebor, \emph{Smoothed analysis of $\kappa(\mathbf{A})$}, J. Complexity, 20 (2004), pp. 97--107.
 

\bibitem{Justin}
P. Zinn-Justin, \emph{Random Hermitian matrices in an external field}, Nucl. Phys. B 497 (1997), pp. 725--732. 
















\end{thebibliography}
\end{document}